\documentclass[11pt]{article}
\usepackage{amssymb}
\newenvironment{pf-of}[1]%
  {\begin{trivlist}\item[\hskip \labelsep{\bf#1}]}%
  {\hfill $\Box$\end{trivlist}}
\def\cal{\mathcal}
\newenvironment{Proof}%
  {\begin{trivlist}\item[\hskip \labelsep{\bf Proof.}]}%
  {\hfill $\Box$\end{trivlist}}
\newcommand{\PSbox}[3]{\mbox{\rule{0in}{#3}\includegraphics{#1}\hspace{#2}}} 
\newtheorem{Theorem}{Theorem}[section]
\newtheorem{Definition}{Definition}[section]

\newtheorem{Lemma}[Theorem]{Lemma}
\newtheorem{Proposition}[Theorem]{Proposition}

\newtheorem{Corollary}[Theorem]{Corollary}
\newtheorem{Conjecture}[Theorem]{Conjecture}
\newtheorem{Claim}[Theorem]{Claim}

\newcommand{\lra}{\longrightarrow}
\def\sqr#1#2{{\vcenter{\vbox{\hrule  height.#2pt
       \hbox{\vrule width.#2pt height#1pt \kern#1pt \vrule width.#2pt}
        \hrule height.#2pt}}}}
\def\bb{\sqr66}   
\let\epsilon=\varepsilon
\let\phi=\varphi
\def\){ \right) }
\def\({ \left( }
\def\[{ \left[ }
\def\]{ \right] }
\def\<{ \langle }
\def\>{ \rangle }
\let\ljunk=\{
\let\rjunk=\}
\def\{{\left\ljunk}
\def\}{\right\rjunk}
\def\p{\partial}

\def\Riem{{\cal R}{\mathrm i}{\mathrm e}{\mathrm m}}

\newcommand{\supp}{{\mathbf S}{\mathrm u}{\mathrm p}{\mathrm p}}

\def\L{\Bbb L}
\def\Vol{{\mathrm V}{\mathrm o}{\mathrm l}}

\def\Tr{{\mathrm T}{\mathrm r}}

\def\ov{\overline}

\newcommand{\cone}{{\mathrm c}{\mathrm o}{\mathrm n}{\mathrm e}}
\newcommand{\R}{{\mathbf R}}

\def\Riem{{\cal R}{\mathrm i}{\mathrm e}{\mathrm m}}

\begin{document} 
\date{ \ } 
\title{Conformal Laplacian and Conical Singularities} \author{Boris
Botvinnik and Serge Preston} \maketitle 
\markboth{{\small Conformal Laplacian and
Conical Singularities}}{{\small B. Botvinnik \&
S. Preston, Conformal Laplacian and
Conical Singularities}}
\vspace{-16mm}

\begin{abstract}
\noindent
We study a behavior of the conformal Laplacian operator $\L_g$ on a
manifold with \emph{tame conical singularities}: when each singularity
is given as a cone over a product of the standard spheres. We study
the spectral properties of the operator $\L_g$ on such manifolds. We
describe the asymptotic of a general solution of the equation $\L_g u
= Q u^{\alpha}$ with $1\leq \alpha\leq \frac{n+2}{n-2}$ near each
singular point.  In particular, we derive the asymptotic of a Yamabe
metric near such singularity.
\end{abstract}
\section{Introduction}\label{int}
{\bf \ref{int}.1. The goal.} The problem we consider in this paper has
two essential parts. Firstly, we study the conformal Laplacian
operator
$$
\L_g= -\Delta_g  + {n-2\over 4(n-1)} R_g 
$$
on a compact Riemannian manifold $(M,g)$ with isolated singularities
of a particular type. Namely, each singular point has a neighborhood
which is a cone over the product of spheres $S^p\times S^q$ endowed
with the standard metric. Secondly, we derive asymptotics for the
positive solutions of the semilinear elliptic equation
\begin{equation}\label{eq_in2}
\L_g u= -\Delta_g u  + {n-2\over 4(n-1)} R_g u = Q u^{\alpha}
\end{equation}
near each singular point. Here $1\leq \alpha\leq
\frac{n+2}{n-2}=\alpha^*$, and $\dim M = n=p+q+1\geq 3$.  We call the
equation (\ref{eq_in2}) the \emph{Yamabe equation}. Indeed, for
$\alpha=\alpha^*$ it corresponds to the \emph{Yamabe problem}.
\vspace{2mm}

\noindent
{\bf \ref{int}.2. Motivation and some prospectives.} Presenting this,
somewhat technical paper, we would like to address and discuss natural
questions which motivate our interest to study the conformal geometry
on manifolds with the cone-type singularities over the product of
spheres.
\vspace{2mm}

First we recall the classical setting for the Yamabe problem. Let
$N$ be a compact smooth manifold, $\dim N = n\geq 3$, and $\Riem(N)$
the space of Riemannian metrics on $N$. We denote by $R_g$ the
scalar curvature and by $d\sigma_{g}$ the volume form for each
Riemannian metric $g\in\Riem(N)$. Then the (normalized)
Einstein-Hilbert functional $I : \Riem(N) \to \R$ is defined as
$$
I : g \mapsto \frac{\int_N R_{g} d\sigma_{g}}
{\Vol_{g}(N)^{\frac{n-2}{n}}}.
$$
Let ${\mathcal C}(N)$ be the space of conformal classes on of
Riemannian metrics on $N$, and $C\in {\mathcal C}(N)$.  The classical
Yamabe problem is to find a metric $\check{g}$ in a given conformal
class $C$ such that the Einstein-Hilbert functional attains its
minimum on $C$: $\displaystyle I(\check{g})=\inf_{g\in C} I(g)$. This
minimizing metric $\check{g}$ is called a \emph{Yamabe metric}, and
the conformal invariant $Y_C(N):=I(\check{g})$ the \emph{Yamabe
constant}.  It is a celebrated result in conformal geometry that the
Yamabe problem has an affirmative solution for closed manifolds, see
\cite{Ya}, \cite{Tr}, \cite{Aubin1}, \cite{Schoen1}.
\vspace{2mm}

The \emph{Yamabe invariant of $N$} is defined as $ Y(N) = \sup_{C\in
{\cal C}(N)} Y_C(N)$. It is well known that the Yamabe invariant
$Y(N)$ is a diffeomorphism invariant of $N$. Furthemore, the Yamabe
invariant completely determines the existence of a metric of positive
scalar curvature: $Y(N)> 0$ if and only if the manifold $N$ admits a
metric of positive scalar curvature. On the other hand, there are just
few examples of manifolds with special properties (in the dimensions
at least than four) for which the value of the Yamabe invariant is
actually known. The fundamental problem here is to compute the Yamabe
invariant $Y(N)$ in terms of other known topological invariants of
$N$.  In particular, it is important to understand a behavior of the
Yamabe constant/invariant under such topological operations as
connected sum of manifolds and, more generally, a surgery.
\vspace{2mm}

O. Kobayashi \cite{Kobayashi} has proven the following estimate for
the Yamabe invariant of a connected sum. Let $N_1$, $N_2$ be two
compact closed manifolds of $\dim N_1=\dim N_2 =n\geq 3$ and
$N_{1,2}=N_1\#N_2$. Then
\begin{equation}\label{surg1}
Y(N_{1,2}) \! \geq \!
\{ \! \begin{array}{cl}
-\(|Y(N_1)|^{n/2}+|Y(N_2)|^{n/2}\)^{2/n}
& \mbox{if}   \ Y(N_1), Y(N_2)\leq 0, 
\\ 
\\
\min\{Y(N_1), Y(N_2)\} &  \mbox{otherwise}.
\end{array}\right.\!\!\!\!\!\!
\end{equation}
Petean and Yun \cite{Petean1} have proven the same formula for the
case when $N_{1,2}=N_1\cup_V N_2$ is a union of two manifolds along a
submanifold $V$ of codimension at least three. This result allowed
Petean to prove that $Y(N)\geq 0$ for any simply connected manifold of
dimension at least five, \cite{Petean2} (see \cite{BR} for the case of
non-simply connected manifolds).
\vspace{2mm}

A difficult case here is to study what is happening with the Yamabe
invariant under surgery if $Y(N)>0$. In more detail, let $S^p\subset
N$ be an embedded sphere with trivial normal bundle. Denote by
$N^{\prime}$ the manifold obtained from $N$ by a surgery along the
sphere $S^p$.  The integer $n-p$ is called a \emph{codimension} of the
surgery. Assume that $Y(N)\leq 0$ and $n-p\geq 3$. Then one can use
(\ref{surg1}) to show that $Y(N^{\prime})\geq Y(N)$, see
\cite{Petean1} and also \cite[Corollary 4.10]{AB}.  In the case when
$Y(N)> 0$ the relationship between the invariants $Y(N)$ and
$Y(N^{\prime})$ is not clear. Consider the trace $W$ of the above
surgery, i.e. $W= N\times I \cup D^{p+1}\times D^{q+1}$, $q=n-p-1$,
where $S^p \times D^{q+1} \subset N\times \{1\}$ is identified with
$S^p \times D^{q+1}\subset \p (D^{p+1} \times D^{q+1})$. In
particular, the boundary of $W$ is a disjoint union of $N$ and
$N^{\prime}$.  According to the elementary Morse theory, one can
choose a Morse function $f: W \to [-1,1]$ with a single nondegenerate
critical point $p$ such that $N=f^{-1}(-1)$, $N^{\prime}=f^{-1}(1)$
and $f(p)=0$. Let $N_{\tau}=f^{-1}(\tau)$, then the manifold
$N_{\tau}$ is diffeomorphic to $N$ if $\tau<0$, and to $N^{\prime}$ if
$\tau>0$. The manifold $N_0$ has an isolated singularity, the vertex
of the cone $C(S^p\times S^q)$ with an appropriate metric, see
Fig. \ref{int}.1.

\PSbox{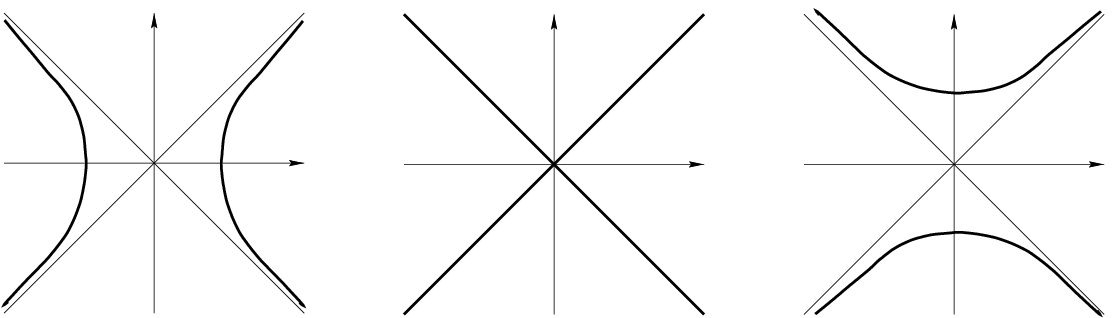}{10cm}{35mm}
\begin{picture}(0,0)
\put(0,10){{\small $N_{\epsilon}$}} \put(-115,10){{\small $N_0$}}
\put(-305,10){{\small $N_{-\epsilon}$}}
\end{picture}

\centerline{{\small {\bf Fig. \ref{int}.1.}}}

\noindent
Hence one can think about a surgery on a manifold $N$ as a deformation
$N_{\tau}$, $-1\leq \tau \leq 1$ through a ``singular point'' $N_0$.
Furthermore, an appropriate metric on the trace $W$ gives a curve
$(N_{\tau},C_{\tau})$ of conformal manifolds. We conjecture that the
function $\tau\mapsto Y_{C_{\tau}}(N_{\tau})$ is a continuous
function. Then under an appropriate choice of conformal classes
$C_{\tau}$, one can approach the following:
\begin{Conjecture}\label{conj1}
Let $N$ be a closed compact manifold of $\dim N=n\geq 5$, and
$N^{\prime}$ be obtained from $N$ by a surgery of codimension at least
three. Then $Y(N^{\prime})\leq Y(N)$.
\end{Conjecture}
If Conjecture \ref{conj1} is indeed true, then the Yamabe invariant
would be a cobordism invariant. In particular, it would mean that if
$N$ is a simply-connected $Spin$-manifold cobordant to zero, then
$Y(N)=Y(S^n)$. 
\vspace{2mm}

From this viewpoint, it is important to understand the behavior of the
conformal Laplacian under surgery, in particular our goal here is to
\emph{extend} the conformal geometry to the category of manifolds with
cone-type singularities over a product of spheres $S^p\times
S^q$. This is the main goal of this paper.
\vspace{2mm}

On the other hand, we think that the results of this paper concerning
the asymptotic of a solution of the Yamabe equation (\ref{eq_in2}) are
of independent value. The asymptotics for singular solutions of the
Yamabe equation were studied thoroughly in the case when the manifold
in question is the standard sphere $S^n$ punctured at $k$ points, see
\cite{MP}, \cite{KMPS}. In particular, this is related to a gluing
construction of metrics of constant scalar curvature under connected
sum operation, see \cite{MPU}. We believe that the results of our
paper could be used to glue metrics of constant scalar curvature under
surgery.
\vspace{2mm}

We have one more application in mind. The paper of the first author
and Akutagawa \cite{AB} gives an affirmative solution of the Yamabe
problem for \emph{cylindrical manifolds}. In particular, in the case
of the \emph{positive cylindrical Yamabe constant} the asymptotic of a
Yamabe metric is \emph{almost conical} near a cylindrical end, see
\cite{AB}.  In that case, it is easy to see that a manifold with
cylindrical ends is equivalent to a manifold with the cone-type
singularities. In particular, the Yamabe problem has a solution for
manifolds with cone-type singularities over the product $S^p\times
S^q$ we study here. As an application, we find an explicit asymptotic
for a Yamabe metric in the specific case related to surgery, see
Corollary \ref{yamabe1}.
\vspace{2mm}

\noindent
{\bf \ref{int}.3. Conformal Laplacian and the Yamabe equation.}
Let $M$ be a compact closed manifold of dimension at least three.  For
a given Riemannian metric $g$ there is the {\it energy functional}
\begin{equation}\label{s-1-1}
E_g(\phi) = \int_M \( |\nabla \phi|^2 + {n-2\over 4(n-1)} R_g \phi^2\)
d\sigma_g.
\end{equation}
Let $1\leq \alpha\leq {n+2\over n-2}=\alpha^*$. We consider the functional
$$
Q_{\alpha}(\phi)= {E_g(\phi) \over \displaystyle \( \int_M
|\phi|^{\alpha+1}d\sigma_g\)^{2\over {\alpha+2}}}, 
\ \ \ \ \phi\in H_2^1(M), \ \ \ \phi>0
$$
with the Euler-Lagrange equation 
\begin{equation}\label{s-1-4}
-\Delta_g u_{\alpha} + {n-2\over 4(n-1)} R_g u_{\alpha} =  
Q_{\alpha}(M) u^{\alpha}_{\alpha}.
\end{equation}
Here $H_2^1(M)$ is
the Sobolev space of functions from $L_{2}(M)$ with their first
derivatives (as distributions) also in $L_{2}(M)$. The Sobolev
embedding theorems imply that for each $\alpha\in [1,\alpha^*)$ there
exists a smooth function $u_{\alpha}>0$ with
$$
\begin{array}{l}
\displaystyle
\int_M u^{\alpha+1}_{\alpha} d\sigma_g =1, \ \ \ 
\mbox{satisfying}
\ \ \ \ 
Q_{\alpha}(u_{\alpha}) = \min\{ Q_{\alpha}(\phi) \ | \ \phi\in H_1(M)\ \}.
\end{array}
$$
We denote $Q_{\alpha}(M) = Q_{\alpha}(u_{\alpha})$. The sign of the
constants $Q_{\alpha}(M)$ is the same for all $1\leq \alpha\leq
\alpha^*$ and it depends only on the conformal class of the metric
$g$, see \cite{Aubin1}.  The existence of a minimizing function in the
case $\alpha=\alpha^*$ is the celebrated Yamabe problem, see
\cite{Ya}, \cite{Tr}, \cite{Aubin0} and \cite{Schoen1}. The solution
of the Yamabe problem had to overcome the fundamental analytic
difficulty concerning Sobolev inequalities for the \emph{critical
exponent} $\alpha^*$.  The constant $Q_{\alpha^*}(M)$ depends only on
the conformal class $[g]$ and is known as the \emph{Yamabe constant}
$Y_{[g]}(M)$ of the conformal class $[g]$.  On the other hand, the
constant $ Q_{1}(M)$ coincides with the first eigenvalue $\mu_1(\L_g)$
of the conformal Laplacian $ \L_g$.
\vspace{2mm}

\noindent
{\bf \ref{int}.4. Manifolds with tame conical singularity.}  Let
$M_0$ be a compact smooth manifold, $\dim M_0 = n\geq 3 $ with the
boundary
$$
\p M_0 = S^p\times S^q \ \ \ \mbox{(where $p+q=n-1$)}.
$$
Let $C(S^p\times S^q)$ be a cone over $S^p\times S^q$ with the vertex
point $x_*$. We glue together $M_0$ and the cone $C(S^p\times S^q)$
along the boundary $S^p\times S^q$ to obtain a manifold with a conical
singularity $x_*\in M$:
$$
M = M_0 \cup_{S^p\times S^q} C(S^p\times S^q).
$$
Now we describe a metric on the manifold $M$.  Let $S^k(r)$ be a
sphere with the standard metric of radius $r$. First, we assume that
$S^p=S^p(r_p)$ and $S^q=S^q(r_q)$. Let $(\theta,\psi,\ell)$ be the
standard coordinate system on the cone $C(S^p\times S^q)$, where
$\theta$, $\psi$ are the spherical coordinates on $S^p$, $S^q$
respectively, and $0\leq \ell\leq \epsilon_1$. In particular,
$(\theta,\psi,0)= x_*$ is the singular point, and
$(\theta,\psi,\epsilon_1)$ give spherical coordinates on $ S^p\times
S^q = \p M_0$. Let $\epsilon_0 < \epsilon_1$. We decompose
$C(S^p\times S^q) = K\cup B$ where
$$
\begin{array}{lcl}
K & = & 
\{ (\theta,\psi,\ell)\in C(S^p\times S^q) \ | 0\leq \ell\leq \epsilon_0\},
\\
\\
B & = &
\{ (\theta,\psi,\ell)\in C(S^p\times S^q) \ | \epsilon_0
\leq \ell\leq \epsilon_1\}.
\end{array}
$$
We denote by $\Lambda = \Lambda(p,q)$ the following constant
$$
\begin{array}{c}
\Lambda=
p(p-1)\frac{2-r_p^{2}}{r_p^{2}}+q(q-1)\frac{2-r_q^{2}}{r_q^{2}}-2pq.
\end{array}
$$ 
We endow $M$ with a Riemannian
metric $g$ satisfying the following properties:
\begin{enumerate}
\item[{\bf (1)}] The scalar curvature function $R_g(x)>0$ if $x\in M_0$.
\item[{\bf (2)}] Let $g_K= g|_K$ be the standard conic metric 
$$
g_K = d\ell^2 + {{{\ell^2}{r_p^{2}}}\over 2} d^2 \theta + 
{{{\ell^2}{r_q^{2}}}\over 2} d^2 \psi
$$
induced from the Euclidian metric in $ C(S^p\times S^q) \subset
\R^{p+1}\times \R^{q+1}, $ where $S^p=S^p(r_p)$, $S^q=S^q(r_q)$. 
(In particular, $\displaystyle
R_{g_K}(\theta,\psi,\ell)=\frac{\Lambda}{\ell^2}$, see Section
\ref{s2}.)
\end{enumerate}
\begin{Definition}\label{def}
{\rm A manifold $(M,g)$ satisfying the above conditions is called a
{\sl manifold with tame conical singularity}.}
\end{Definition} 
The conformal Laplacian $\L_g = - \Delta_g + {n-2\over 4(n-1)}R_g$ is
well-defined on the manifold $M$ without the singular point $x_*$. In
particular, we have the Yamabe equation $\L_gu=Q_{\alpha}u^{\alpha}$
on $M\setminus \{x^*\}$ for $1\leq \alpha \leq \alpha^*={n+2\over
n-2}$.
\vspace{2mm}

\noindent
{\bf \ref{int}.5. Results.} We study the following issues:
\begin{enumerate}
\item[{\bf (A)}] The spectral properties of the conformal Laplacian $\L_g$ on
a manifold with tame conical singularity.  
\item[{\bf (B)}] The asymptotic of a general solution of the linear
Yamabe equation (i.e. when $\alpha=1$) near the singular point.
\item[{\bf (C)}] The asymptotic of a general solution of the non-linear
Yamabe equation near the singular point for $1< \alpha \leq \alpha
^{*}$.
\end{enumerate} 
The paper is organized as follows. We describe in detail the geometry of
$M$ and the Yamabe equation near the singular point in Section
\ref{s2}. We define appropriate weighted Sobolev spaces in Section
\ref{s3}. We study the conformal Laplacian $\L_g$ on $M$ in Section
\ref{s4}, in particular, we prove that under some dimensional
restrictions the operator $\L_g$ is positive. We analyze the 
functional $I_{\alpha}$ on $M$ and prove a weak version of the Yamabe
theorem in Section \ref{s5}. We study the asymptotic of a general
solution of the linear and nonlinear Yamabe equations in Sections
\ref{s6} and \ref{s7} respectively. We put together some technical
results and calculations in Appendix (Section \ref{ap}). 
\vspace{2mm}

\noindent
{\bf \ref{int}.6. Acknowledgments.} The authors are grateful to Kazuo
Akutagawa for useful discussions on the Yamabe invariant and to the
referee for his helpful and encouraging remarks. We are also grateful
to Ms. Melisa Gilbert for the editorial help.
\section{The Yamabe equation on the cone}\label{s2}
The cone $K \subset \R^{p+1}\times \R^{q+1}=\R^{n+1}$ is given by
\begin{equation}\label{cone-1}
\!\!\!\begin{array}{l}
\!\!\!K\! = \!\{\! (x_1,\ldots,x_{n+1}) |
\frac{1}{r_q^2}\left(x_1^2+\cdots+x_{p+1}^2\right)\!- \!\frac{1}{r_q^2}
\left(x_{p+2}^2+\cdots+x_{n+1}^2\right)\! =\!0 \!\}. 
\end{array}\!\!\!\!\!\!\!\!
\end{equation}
We use the polar coordinates $(r,\theta)$ in $\R^{p+1}$, where $\theta\in
S^p$, and $(\rho,\psi)$ in $\R^{q+1}$, where $\psi\in S^q$. Then the
equation (\ref{cone-1}) can be written as $\frac{\rho}{r_p}
=\frac{r}{r_q}$.
\vspace{6mm}

\noindent
\parbox{2.4in}{
\hspace*{1mm}\PSbox{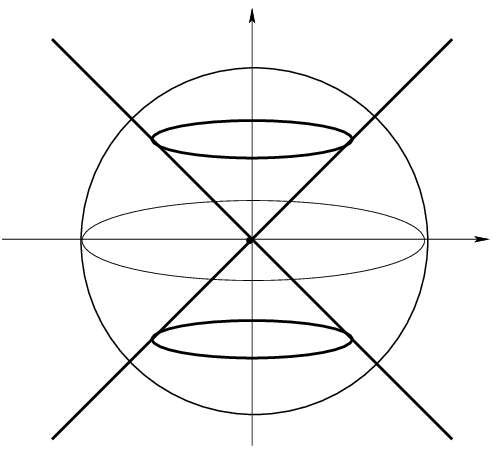}{10cm}{40mm}
\begin{picture}(0,0)
\put(125,115){{\small $K$}}
\put(80,135){{\small $x^{p+2},\!\ldots,\!x^{p+q+2}$}}
\put(110,85){{\small $x^{1},\!\ldots,\!x^{p+1}$}}
\end{picture}
\centerline{{\small {\bf Fig. \ref{s2}.1.}}}}
\parbox{2.5in}{
The coordinate $\ell$ is the distance from a point on the cone to the
vertex $x^*$. Then $\ell=r\sqrt{2}=\rho\sqrt{2}$. The embedding $K
\subset
\R^{n+1}$ induces the metric $g_K$ from the standard metric on
$\R^{n+1}$.  In polar coordinates the metric $g_K$ is given as
$$
\begin{array}{lcl}
g_K &=& \displaystyle
d r^2 + r^2   d^2 \theta  + r^2   d^2 \psi
\\
\\
&=& \displaystyle
d\ell^2 + \frac{r_p^2\ell^2}{2}d^2 \theta + 
\frac{r_q^2\ell^2}{2}d^2 \psi.
\end{array}
$$
}

\noindent
We compute the determinant $|g_K|$:
$$
\begin{array}{c}
|g_K| = \displaystyle
\det(g_{ij})= 1\cdot \( \frac{r_p^2\ell^2}{2}\)^p\cdot |g_{\theta}|
\cdot \(\frac{r_q^2\ell^2}{2}\)^q\cdot |g_{\psi}|, \ \ \mbox{and}
\\
\\
\sqrt{|g_K|} = \displaystyle
\frac{r_p^p r_q^q\ell^{n-1}}{2^{\frac{n-1}{2}}} \cdot \displaystyle
\sqrt{|g_{\theta_{\mbox{\ }}}^{\mbox{\ }}|}\displaystyle
\cdot \sqrt{|g_{\psi}|} .
\end{array}
$$
Here $g_{\theta}$ and $g_{\psi}$ are the standard metrics on the
spheres $S^p(r_p)$ and $S^q(r_q)$ respectively. In particular, we have
the following formula for the volume element:
\begin{equation}\label{cone-4}
d\sigma_{g_K} = {r_p^{p}r_q^{q}\ell^{n-1}\over 2^{{n-1\over 2}}} d\ell
\wedge d\sigma_{g_{\theta}} \wedge d\sigma_{g_{\psi}}.
\end{equation}
The Laplace operator on the cone is then given as
$$
\Delta_{g_K} = 
{1\over \sqrt{|g_K|}}\p_{x^i}\(\sqrt{|g_K|}\ g^{ij}_K \p_{x^j}\).
$$
The matrix $\(g^{ij}_K\)$ is block diagonal:
$$
\(g^{ij}_K\) =
\(\mbox{\begin{picture}(80,50)
\put(15,0){\line(0,1){50}}
\put(0,30){\line(1,0){50}}
\put(15,0){\line(1,0){65}}
\put(50,30){\line(0,-1){65}}
\put(5,35){$1$}
\put(25,10){$\frac{2g^{ij}_{\theta}}{\ell ^{2}r_p^{2}}$}
\put(60,-20){$\frac{2g^{ij}_{\psi}}{\ell ^{2}r_q^{2}}$}
\end{picture}}\) \ .
$$
We obtain:
$$
\begin{array}{rcl}
\Delta_{g_K} &=& \displaystyle
{2^{{n-1\over 2}}\over \ell^{n-1}}\p_{\ell}\( {\ell^{n-1}\over
2^{{n-1\over 2}}}\p_{\ell}\) + {2\over {r^{2}_p\ell^2}}\Delta_{\theta}
+ {2\over {r^{2}_q\ell^2}}\Delta_{\psi}
\\
\\
&=& \displaystyle
{\p^2\over \p\ell^2} + (n-1){1\over \ell}{\p \over \p \ell}
+ {2\over r^{2}_p\ell^2}\Delta_{\theta} + {2\over r^{2}_q\ell^2}\Delta_{\psi} .
\end{array}
$$
The metric $g_K$ is a particular case of a ``double wrapped metric''
on the product $I\times S^p\times S^q$. In general, a double wrapped metric
is given as
$$
d\ell^2 + \phi_p^2(\ell) d^2 \theta + \phi_q^2(\ell) d^2 \psi.
$$
In our case 
$$
\phi_p(\ell)={r_p\ell\over \sqrt{2}}, \ \ \phi_q(\ell)={r_q\ell\over \sqrt{2}},
\ \
 \dot{\phi}_p = {r_p\over \sqrt{2}}, \ \ \dot{\phi}_q={r_q\over \sqrt{2}},
\ \ 
\ddot{\phi}_p= 0 , \ \ \ddot{\phi}_q= 0 .
$$
We choose the orthonormal bases for the tangent space to $K$:
$$
F_0=\p_{\ell}, \ \ \ F_i, \ i\in\overline{p}, \ \ \
F_j, \ j\in\overline{q},
$$
where $\overline{p}=\{2,\ldots,p\}$, $\overline{q}=\{p+1,\ldots,p+q\}$
are orthonormal basis for the standard metrics on $S^{p}$ and $S^{q}$
correspondingly.  Then we have the following formulas for the
curvature operator ${\mathcal R}$ acting on the bundle of 2-forms
(see, for instance the book by P. Petersen
\cite{Petersen}):
$$
\begin{array}{lcl}
{\mathcal R}(F_0\wedge F_i) &=& 0 , \ \ \ i\in \overline{p},
\\
\\
{\mathcal R}(F_0\wedge F_j) &=& 0 , \ \ \ i\in \overline{q},
\\
\\
{\mathcal R}(F_{i_1}\wedge F_{i_2}) &=& \displaystyle
{{2-r_p^{2}}\over {r_p^{2}\ell^2}} F_{i_1}\wedge F_{i_2},
\ \ \ i_1,i_2\in \overline{p},
\\
\\
{\mathcal R}(F_{j_1}\wedge F_{j_2}) &=& \displaystyle
{{2-r_q^{2}}\over {r_q^{2}\ell^2}} F_{j_1}\wedge F_{j_2},
\ \ \ j_1,j_2\in \overline{q},
\\
\\
{\mathcal R}(F_{i}\wedge F_{j}) &=& \displaystyle
- {1\over \ell^2 } F_{i}\wedge F_{j},
\ \ \ i\in \overline{p}, \ \ j\in \overline{q} .
\end{array}
$$
Thus we have:
\begin{equation}\label{cone-6}
\begin{array}{lcl}
{\mathcal R}ic(F_0) &=& 0,
\\
\\
{\mathcal R}ic(F_i) &=& \displaystyle
\((p-1){{2-r_p^{2}}\over {r_p^{2}\ell^2}} - q{1\over \ell^2}\)F_i ,
\\
\\
{\mathcal R}ic(F_j) &=& \displaystyle
\((q-1){{2-r_q^{2}}\over {r_q^{2}\ell^2}} - p{1\over \ell^2}\)F_i .
\end{array}
\end{equation}
We compute the scalar curvature:
$$
\begin{array}{rcl}
R_{g_K}& = & \displaystyle
\Tr\ {\mathcal R}ic = p\left[ (p-1)\frac{2-r_q^{2}}{r_p^{2}\ell ^{2}}-
\frac{q}{\ell^{2}}\right] + 
q\left[ (q-1)\frac{2-r_q^{2}}{r_q^{2}\ell ^{2}}-\frac{p}{\ell^{2}}\right] =
\\
\\
& = & \displaystyle {1\over \ell^2}\cdot \left[
p(p-1)\frac{2-r_p^{2}}{r_p^{2}}+ q(q-1)\frac{2-r_q^{2}}{r_q^{2}}-2pq
\right] = {\Lambda\over\ell^2}, \ \ \mbox{with}
\\
\\
\Lambda&=&\displaystyle p(p-1)\frac{2-r_p^{2}}{r_p^{2}}+
q(q-1)\frac{2-r_q^{2}}{r_q^{2}}-2pq .
\end{array}
$$
We refer to the Appendix for some properties of $\Lambda $ as a function
of $p$, $q$, $r_p$, $r_q$.
\vspace{2mm}

Now we can rewrite the Yamabe equation $\L_g u = Q_{\alpha}u^{\alpha}$
for the cone $K$. We obtain:
$$
\begin{array}{c}
\displaystyle
{\p^2 u\over \p  \ell^2} + {n-1\over \ell}{\p u \over \p\ell}+
{2\over {r_p^{2}\ell^2}}\Delta_{\theta}u +
{2\over {r_q^{2}\ell^2}}\Delta_{\psi}u
- {n-2\over 4(n-1)} {\Lambda\over \ell^2} u +
Q_{\alpha}u^{\alpha}=0 ,
\\
\\
\displaystyle
\mbox{where $1\leq \alpha\leq
\alpha^*$ and with}
\ \ \ \ \int_K u^{\alpha+1} <1, \ \ \
u_{\alpha}>0.
\end{array}
$$ 
We notice that $\int_K u^{\alpha+1} <1$ since
$\int_M u^{\alpha+1} =1$.
\section{Sobolev spaces on $M$}\label{s3}
In this section we introduce appropriately weighted Sobolev spaces on a
manifold with tame conical singularities and review their basic
properties.  We denote $d\sigma_g$ the volume form on $M$
corresponding to the metric $g$. According to the computations above,
the form $d\sigma_g$, restricted on the cone part $K$, is given by
$$
d\sigma_g = {\ell^{n-1}\over 2^{{n-2\over 2}}}
\sqrt{|g_{\theta}^{\mbox{ \ }}|} \sqrt{|g_{\psi}|}d\ell\wedge
d\sigma_{\theta} \wedge d\sigma_{\psi}.
$$
Let $L_2(M)$ be the Hilbert space of functions $\phi$ on $M$.
Clearly, the restriction $\phi|_K=\phi(\ell,\theta,\psi)$ of a function
$\phi\in L_2(M)$ on the cone $K$ satisfies the property
\begin{equation}\label{sob-1} \int_{0}^{\epsilon_0}
\int_{S^p}\int_{S^q} |\phi(\ell,\theta,\psi)|^2 \ell^{n-1} d\ell
\wedge d\sigma_{\theta} \wedge d\sigma_{\psi} < \infty.
\end{equation}
The proof of the following lemma is standard.
\begin{Lemma}\label{sob-1b}
Let $\phi$ be a function on $M$ with $\int_{M\setminus K} |\phi|^2
d\sigma_g<\infty$, and with the asymptotic behavior $\phi|_K=
\phi(\ell,\theta,\psi)=\ell^s(1+O(1))$ as $\ell\rightarrow 0$.
Then the following statements are equivalent:
\begin{enumerate}
\item[{\bf (1)}] The integral {\rm (\ref{sob-1})} converges.
\item[{\bf (2)}] The function $\phi\in L_2(M)$.
\item[{\bf (3)}] $s> -{n\over 2}$.
\end{enumerate}
\end{Lemma}
Let $\chi \in C^{\infty}(M\setminus x_*)$ be a positive weight
function satisfying
$$
\chi(x) = \{\begin{array}{ll}
1 & \mbox{if $x\in M_0$},
\\
{1\over \ell} & \mbox{if $x\in K$}
\end{array}
\right. \ \ \ \mbox{and} \ \ \ 0<\chi(x)\leq1 \ \ 
\mbox{for all $x\in M \setminus x_*$}.
$$
It is easy to construct such a function. Then a \emph{weighted Sobolev
space} $H^k_2(M)=H^k_2(M,g)$ is defined as a space of locally
integrable functions $\phi\in L^{loc}_1(M,g)$ with finite norm
\begin{equation}\label{sob-1a}
\|\phi\|^2_{H_k^2(M)} = 
\int_M\(\sum_{i= 0}^{k}\chi^{2(k-i)}\sum_{|\mu|= i} |D^{\mu}\phi|^2\)
d\sigma_g.
\end{equation}
The derivatives here are understood as distributions (i.e.
$D^{\mu}\phi\in L^{loc}_1(M,g)$, $|\mu|\leq k$).  We denote by
$C^{\infty}_*(M)$ the space of smooth functions on $M$ that are equal
to zero in some neighborhood of the singular point $x_*$.  We notice that
although the norm of a function $\phi$ in Sobolev spaces defined here
depends on the choice of the weight $\chi $, the property of $\phi$ 
belonging to the corresponding Sobolev space does not.
\vspace{2mm}

\noindent
{\bf Remark.} The spaces $H^k_2(M,g)$ defined above coincide with the
spaces $W^k_2(M,1,\rho^{2k})$ defined by H. Triebel \cite[Section
3.24, 3.26]{Triebel}, where $\rho(x)= \chi(x)$.
\vspace{2mm}

\noindent
In particular, according to \cite[Section 3.2.3]{Triebel}, the
spaces $H^k_2(M,g)$ are complete Hilbert spaces with the scalar
product
$$
\<\phi,\psi\>_{H_k^2(M)} = 
\int_M\(\sum_{i= 0}^{k}\chi^{2(k-i)}\sum_{|\mu|= i} D^{\mu}\phi
D_{\mu}\psi\) d\sigma_g . 
$$
The spaces $H^k_2(M,g)$ are closely related to the
functional spaces $V^k_{2,0}(M,g)$ and $W^k_{2,\beta}(M)$ defined in
the book
\cite{KMR}.  The space $V^k_{2,0}(M,g)$ is a closure of the space
$C^{\infty}_*(M)$ in the norm (\ref{sob-1}), \cite[Section
6.1.1]{KMR}, and the space $W^k_{2,\beta}(M,g)$ is defined as a space
of locally integrable functions $\phi\in L_{loc}^1(M,d\sigma_g)$ with
the following finite norm:
$$
\|\phi\|^2_{W^{\ell}_{2,\beta}(M)} = 
\(\int_M r^{2\beta}\sum_{|\mu|\leq \ell} |D^{\mu}\phi|^2 
d\sigma_g\)^{1/2},
$$
where the weight function $r$ is related to our weight function as
$\chi\sim r^{-1}$ (see \cite[Section 7.1.2]{KMR}). Clearly we have the
following embeddings
\begin{equation}\label{sob-4}
V^k_{2,0}(M,g)\subset H_k^2(M,g) \subset W^k_{2,0}(M,g) .
\end{equation}
These embeddings are continuous; moreover, if $\beta> k- n/2$, the
spaces $V^k_{2,\beta}(M,g)$ and $W^k_{2,\beta}(M,g)$ coincide and
their norms are equivalent (see \cite[Theorem 7.1.1]{KMR}).  In our
case, when $\beta= 0$ (when $k<n/2$) three spaces in (\ref{sob-4})
coincide. This implies the following result.
\begin{Lemma}\label{sob-5}
Let $n= \dim M\geq 5$. Then the space $C^{\infty}_*(M)$ is dense in the 
spaces
$H^1_2(M,g)$ and $H^2_2(M,g)$.
\end{Lemma}
{\bf Remark.} Lemma \ref{sob-5} is the first point when the analysis
implies the dimensional restriction $n\geq 5$. We do not know whether
some modification of Lemma \ref{sob-5} holds in dimensions $n=3,4$.
\vspace{2mm}

Now we recall the embedding theorem, following the exposition of
H. Triebel \cite{Triebel}, for the spaces $H^k_2(M,g)=  W_2^k(\Omega,
\rho^0,\rho^2)$ (where $\rho= \chi$).
\begin{Theorem}\label{sob-6}
Let $M$ be a manifold with a tame conical singularity as above, $n\geq
5$.  Let $2\leq p \leq 1 + \alpha^*$, $\alpha^*= {n+2\over n-2}$. Then
there exists a continuous linear embedding
$$
i_p: H^k_2(M,g) \subset L_p (M)
$$
for $k\geq 1$. Furthermore, the embedding operator $i_p$ is compact for
$2\leq p < 1 + \alpha^*$.
\end{Theorem}
\begin{Proof}
This result follows by combining the embedding theorems of Triebel
(see the book \cite[Theorems 3.5.1, and 3.8.3]{Triebel}) for the
domains in $\R^n$ with conical points on the boundary with the
classical embedding theorems for Sobolev spaces on compact manifolds
(see, for example, \cite[Theorems 2.10, 2.33]{Aubin1}).
\end{Proof}
There is an important special case here. Clearly, we have the
embeddings:
$$
H_2^2(M,g) \subset H_2^1(M,g) \subset L_2(M).
$$
Theorem \ref{sob-6} implies that the embedding $H_1^1(M,g) \subset
L_2(M)$ is compact. We will need the following lemma which follows
from Theorem \ref{sob-6}.
\begin{Lemma}\label{sob-7}
Let $n\geq 5$, then in the space $H_2^1(M,g)$ the norms
$$
\|\nabla \phi\|_2^2 + \|\chi\phi\|_2^2 \ \ \ \mbox{and} \ \ \
\|\nabla \phi\|_2^2 + \|\phi\|_2^2
$$
are equivalent, in particular, for some $C>0$
$$
\|\chi\phi\|_2^2 \leq C\(\|\nabla \phi\|_2^2 + \|\phi\|_2^2\).
$$
\end{Lemma}
The following lemma allows us to compare norms in $L_2(M_0)$ and
$L_2(M)$.
\begin{Lemma}\label{new}
Let $n\geq 5$. Then for any constant $a>0$ there exists a constant
$C>0$ such that
$$
\|\nabla u\|_{L_2(M)} + a \| u\|_{L_2(M_0)} \geq C \| u\|_{L_2(M)}
$$
for any function $u\in H^1_2(M)$.
\end{Lemma}
\begin{Proof}
Notice, first, that for a function $u$ as in the formulation of Lemma
\ref{new}, we have $u|_{M_0}\in L_2(M_0)$.  Assume that for some $a>0$
there is no such constant $C>0$. Then there exists a sequence $u_n\in
H^1_2(M)$, such that
$$
\|u_n\|_{L_2(M)} =  1, \ \ \mbox{and} \ \ \ \|\nabla u_n\|_{L_2(M)}
\to 0,
\ \ \ \| u_n\|_{L_2(M_0)}  \to 0
$$
as $n\rightarrow \infty$. Recall that the embedding $H^1_2(M)\subset
L_2(M)$ is compact; thus, passing to a subsequence, if necessary, we
get $u_n \rightarrow u_{\infty}$ in $L_2(M)$ for some $u_{\infty}\in
L_2(M)$. Clearly we have $\|u_{\infty}\|_ {L_2(M)} =  1$, but
$u_{\infty}|_{M_0}= 0$ almost everywhere since $u_n\rightarrow 0$ in
$L_{2}(M_{0})$ as
$n\rightarrow \infty$. Now we choose a function $\phi\in
C^{\infty}_0(M)$ (in particular, $\phi$ has compact support). We have
$$
|\<\nabla u_n,\nabla\phi\>_{L_2(M)} | \leq \| \nabla u_n\|_{L_2(M)}\cdot
\| \nabla \phi \|_{L_2(M)}.
$$
As $n\rightarrow \infty$ the right-hand side goes to zero, thus the
left-hand side too. However
$$
\<\nabla u_n,\nabla\phi\>_{L_2(M)} =  \< u_n,\Delta \phi\>_{L_2(M)}
\lra
\< u_{\infty},\Delta \phi\>_{L_2(M)}
$$
(using that $\supp(\phi)$ is compact).  Thus, $\<u_{\infty },\Delta \phi
\>_{L_{2}(M)}=0$ for all $\phi \in C^{\infty }_{0}(M)$, i.e.
$u_{\infty}\in L_2(M)$ is a weak solution of the equation $\Delta_{g}
u = 0$. Since the Laplacian $\Delta_{g}$ is an elliptic operator,
$u_{\infty}$ is analytic in $M$. Since $u_{\infty}|_{M_0}\equiv 0$,
then $u_{\infty}\equiv 0$ on $M$. This contradicts to 
$\Vert u_{\infty }\Vert _{L_{2}(M)}=1$.
\end{Proof}
Let $H^{k}_{2,*}(M,g)$ be the corresponding Sobolev space of functions
on $M$ with some open neighborhood of the singular point  $x_*$ removed.
\begin{Lemma}\label{new-l1}
Let $f=f(l,\theta,
\psi) \in H^{k}_{2,*}(M) $ be a measurable function,  so that 
$f\in H^{k}_{2}(M_*)$, and $f$ has asymptotic behavior $f \sim l^{q}$ as $ l\longrightarrow
0$ near the point $x_*$. Then, $f\in H_{2}^{k}(M)$ if and only if
$k<q+\frac{n}{2}$.
\end{Lemma}
\begin{Proof} We have that the integral in the formula 
(\ref{sob-1a}) near $\ell= 0$ has the form
$$
\int _{0}^{\epsilon} \ell^{2(q-k)}\ell^{n-1}d\ell.
$$
The integral in (\ref{sob-1a}) (taken over $M\setminus U$ where $U$ is
some open neighborhood of $x_*$) is finite by the condition on $f$.
It is easy to see that the Sobolev norm of the function $f$ over all
$M$ in (\ref{sob-1a}) is finite if and only if $2(q-k)+n-1 >-1$, or,
if $q-k+\frac{n}{2}>0$. This gives the condition $k< q+\frac{n}{2}$.
\end{Proof}
Later it will be convenient for us to refer to a class of functions
with particular asymptotic near the singular point $x_*$.
\begin{Definition}\label{new-d1}
{\rm Let $1 \leq \alpha \leq \alpha ^{*}$. A locally integrable
function $f\in L_{1,*}(M)$ is called an {\sl $\alpha$-basic function}
if $f\sim \ell^{-\frac{2}{\alpha -1}}$ as $ \ell\to 0$,
namely, 
$$
c \cdot \ell ^{-\frac{2}{\alpha -1}}\leq  
\vert f(\ell , \theta, \psi
)\vert \leq C \cdot \ell ^{-\frac{2}{\alpha -1}}
$$ 
for some positive constants $c$, $C$.}
\end{Definition}
In particular, if $\alpha = \alpha ^{*}= \frac{n+2}{n-2}$, an
$\alpha$-basic function means that $f \sim \ell^{-\frac{n-2}{2}}$ near
$x_*$.
\vspace{2mm}

It is easy to see that Lemma \ref{new-l1} implies the following
result:
\begin{Proposition}\label{new-p1}
Let $1 \leq \alpha \leq \alpha ^{*}$. Then an $\alpha$-basic function
$f$ belongs to $L_{2}(M)$ if and only if $\frac{n+4}{n} < \alpha \leq
\alpha ^{*}$. Furthermore, $f$ does not belong to $H^{1}_{2}(M)$ if
$\alpha\in [1,\alpha ^{*})$.
\end{Proposition}
\begin{Proof} Let $f$ be an $\alpha$-basic function. Then 
$f\in H_{2}^{k}(M)$ if and only if
 $k < \frac{n}{2}-\frac{2}{\alpha -1}$. Then 
$k= 0$ gives  $\alpha -1 >\frac{4}{n}$ or $\alpha
>\frac{n+4}{n}$, and   $k= 1$ gives $\frac{2}{\alpha -1} <
\frac{n}{2}-1= \frac{n-2}{2}$, or $\alpha -1 >\frac{4}{n-2}$. But this
is exactly opposite to the main condition for the parameter
$\alpha$: $\alpha \leq \alpha ^{*}$.
\end{Proof}
\section{Spectral properties of the conformal Laplacian}\label{s4} 
{\bf \ref{s4}.1. General remarks.} In this section we study 
the conformal Laplacian  
\begin{equation}\label{n1}  
\L_g = -\Delta_g+ c(x), \ \ \ \mbox{with} \ \ \ c(x)= {n-2\over
4(n-1)}R_g(x)
\end{equation}  
on a manifold $M$ with tame conical singularity.  
\vspace{2mm}

The classical Laplace operator $\Delta_g$ has been studied thoroughly
on compact manifolds and on manifolds with certain singularities.  In
particular, the spectral properties of the Laplacian and its heat
kernel are well-understood.  The conformal Laplacian is studied
predominantly through its relation to the conformal geometry in
general, and the Yamabe problem in particular, although there has been
important work done on its general properties as well (for example,
work by T. Parker, S.  Rosenberg \cite{Parker-Rosenberg}, by
T. Branson \cite{Branson}, see also the references given in the book
by P. Gilkey \cite{Gilkey}).
\vspace{2mm}

On the other hand, operators of the type (\ref{n1}) were studied much
as \emph{Schr\"odinger} \emph{operators}. For example, the book by S.
Mizohata \cite[Chapter 8]{Mizohata}) studies the Schr\"odinger
operator $-\Delta+ c(x)$ on $\R^3$ where the function $c(x)$ has
a singularity at the origin of order $\displaystyle |x|^{-{3\over
2}+\epsilon}$ or weaker.  Under these conditions, the operator
$-\Delta+ c(x)$ is bounded from below and has a unique self-adjoint
extension (\emph{Friedrichs' extension}). Additional asymptotic
conditions (when $|x| \to \infty$) on $c(x)$ ensure that the spectrum
of $-\Delta+ c(x)$ is discrete and of finite multiplicity.
\vspace{2mm} 

\noindent 
{\bf \ref{s4}.2. Basic results.}  We start with the following standard
property of $\L_g$ (this follows from the basic results on the Sobolev
spaces).
\begin{Proposition}\label{Pn1} The conformal Laplacian $\L_g=
-\Delta_g + {n-2\over 4(n-1)} R_g$ is densely defined in $\L_{2}(M)$
with the domain $H_{2}^{2}(M,g)$. Furthermore, the operator $\L_g$ is
symmetrical and continuous with respect to the norm in
$H_{2}^{2}(M,g)$.
\end{Proposition} 
\begin{Proof} It is obvious that
this operator is symmetric with the domain $C^{\infty }_{*}(M)$. We
look at the norm 
$$ 
|\L_g\phi|^2 = \(-\Delta_g\phi + {n-2\over 4(n-1)}R_g\phi\)^2 \leq 2
\(\(\Delta_g\phi\)^2 + \({n-2\over 4(n-1)}\)^2 R_g^2 \phi^2\).
$$ 
Integration over $M$ gives: 
$$
\|\L_g\phi\|^2_{L_2(M)} \leq 2\( \int_M \(\Delta_g\phi\)^2 d\sigma_g +
C \int_M R_g^2 \phi^2 d\sigma_g\).  
$$ 
Notice now that $|R_g(x)|\leq C_1\chi(x)$ for some constant $C_1$
(Indeed, on $M\setminus K$ the scalar curvature is bounded from above
while on the cone $K$ it is equal to $\Lambda\chi(x)$.)  Thus using
the definition of the norm (\ref{sob-1a}) we obtain 
$$
\|\L_g\phi\|^2_{L_2(M)} \leq 2\( \int_M \(\Delta_g\phi\)^2d\sigma_g +
C_1 \int_M \chi^2 \phi^2d\sigma_g\) \leq C_g \|\phi\|^2_{H_2^2(M,g)}.
$$ 
This proves Proposition \ref{Pn1}.
\end{Proof} 
\begin{Theorem}\label{Th1new} 
Let $(M,g)$ be a manifold with tame
conical singularity of $\dim M= n \geq 3 $. Then the quadratic form
$\<\L_g u,u \>$ defined on $H_2^1(M,g)$ is such that 
$$ 
D \|u
\|_{H^1_2(M,g)}^2 \leq \<\L_g u, u \> \leq C \|u \|_{H^1_2(M,g)}^2
$$  
for some constants $C>0$ and $D>-\infty$ (i.e. the form $\<\L_g u,u \>$
is bounded from below).  Furthermore,
\begin{enumerate}  
\item[{\bf (1)}] if $\Lambda>0$, then  
\begin{equation}\label{pos}  C_1\|u 
\|_{H^1_2(M,g)}^2 \leq \<\L_g u, u \> \leq  C_2 \|u 
\|_{H^1_2(M,g)}^2  
\end{equation}  
for some constants $C_1>0$, $C_2>0$;
\item[{\bf (2)}] if $\Lambda \leq 0$ then for each $r_p$, 
$r_q$  there is only a finite number of nonpositive eigenvalues (each of 
finite  multiplicity) of the operator 
$\L_g$.  
\end{enumerate}  
\end{Theorem}  
\begin{Corollary}\label{Cn3} 
Let $\Lambda>0$. Then for any $u\in H_2^1(M,g)$  the operator $\L_g$ 
satisfies  $$  \|\L_g(u)\|_{L_2(M)} \geq C \|u\|_{L_2(M)}  $$  for 
some constant $C>0$. 
\end{Corollary} 
\begin{Corollary}\label{Cn4}  
There exists the 
Friedrichs' self-adjoint extension $\widetilde{\L}_g$  of the operator 
$\L_g$. The extension $\widetilde{\L}_g$ is  semi-bounded with the same 
bounding constant. The range of  $\widetilde{\L}_g$ coincides with 
$L_2(M)$.  
\end{Corollary}  
\begin{Proof}  This follows directly from 
4.2 and the Neumann Theorem (see \cite[Theorem 17]{Egorov}).
\end{Proof}  
\begin{Corollary}\label{Cn5}  
The conformal Laplacian $\L_g$ is essentially self-adjoint, and its
self-adjoint extension is unique.
\end{Corollary}  
\begin{Proof}  It follows from the fact 
that $\L_g$ is strictly positive and symmetric, see \cite[Theorems 28,
29]{Egorov}.  
\end{Proof} 
{\bf Remark.} We have that the range $R(\widetilde{\L}_g)=
L_2(M)$. Thus, the range $R(\L_g)$ is dense in $L_2(M)$. Thus,
essential self-adjointness of $\L_g$ follows also from
\cite[Lemma 8.14]{Mizohata}. A direct proof of the density of
$R(\L_g)$ in $L_2(M)$ may be given following the proof (presented in
\cite[Chapter 8, Section 13]{Mizohata}) for the Schr\"odinger operator
$-\Delta + c(x)$. The latter one utilizes the asymptotic estimates of
the Green function for $\L_g$ similar to one obtained by V. Maz'ya,
S. Nasarow, B. Plamenewski (see \cite{Maz}).
\begin{Theorem}\label{Cn6} 
Let $\Lambda >0$. Then the self-adjoint extension $\widetilde{\L}_g$
of the conformal Laplacian $\L_g$ has discrete positive spectrum of
finite multiplicity.
\end{Theorem} 
\begin{Proof} We use compactness of the embedding 
${\cal D}(\widetilde{\L}_g)= H_2^2(M,g)\subset L_2(M)$ and Rollich
Theorem (see, say, \cite[Theorem 3.3.]{Mizohata}) to prove that the
inverse operator $\widetilde{\L}_g^{-1}$ is compact and has discrete
spectrum $\{\mu_n\}$ ($\lambda_n^{-1}= \mu_n\rightarrow 0$ as
$n\rightarrow \infty$) of finite multiplicity. Also $\<\L_gu,u\>\geq
C\|u\|_{L_2(M)}$ gives that $\sigma(\widetilde{\L}_g)\subset \R_+$.
\end{Proof} 
{\bf Remark.}  One can also prove that ${\cal D}(\widetilde{\L}_g)=
H_2^2(M,g)$ using arguments similar to those used by H. Triebel
\cite[Theorem 6.4.1]{Triebel} for the operator $-\Delta +\chi^2$.
\begin{pf-of}{\ref{s4}.3. Proof of Theorem \ref{Th1new}.}  The upper
estimate follows from the symmetry of the form $\<\L_{g}u,u\>$ and the
definition of the Sobolev spaces.  
\vspace{2mm}

For the proof of the lower estimate, consider first an arbitrary
smooth function $u\in C_*^{\infty}(M)$. We have
\begin{equation}\label{c-1}
\begin{array}{rcl} 
\!\! \<\L_g u, u \> &= & \displaystyle 
\int_M\!\(\!- u \cdot
\Delta_g u\! +\! {n-2\over 4(n-1)} R_g u^2 \!\) d\sigma_g 
\\ \\ &= &
\displaystyle 
\int_M\!\( |\nabla_g u|^2\! + \!{n-2\over 4(n-1)} R_g u^2\! \)
d\sigma_g.  
\end{array} 
\end{equation} 
Here we use $u\in C_*^{\infty}(M)$, so $u$ is zero in some
neighborhood of $x_*$.  Indeed, let $u\equiv 0$ in $\{|\ell|,
\nu\}\subset K \subset M$ for some $\nu$. Then $M \setminus
\{|\ell|<{\nu\over 2}\}$ is a manifold with the boundary $\{|\ell|=
{\nu\over 2}\}$, and $u$ is zero in a neigborhood of this boundary.
\vspace{2mm}

We decompose $M$ as follows: 
$$
\begin{array}{l}\displaystyle M = \( M_0 \cup B_{\sigma}\) \cup
\((B\setminus B_{\sigma})\cup K\), \ \ \mbox{with} \\ \\ \displaystyle
B_{\sigma} = \{ x= (r,\ell) \in C(S^p\times S^q) \ | \ \sigma - r \ell
|\ell| <\sigma \}, 
\end{array} 
$$ 
so that the scalar curvature $R_g\geq R_0>0$ on $M_0 \cup B_{\sigma}$,
while on $(B\setminus B_{\sigma})\cup K$ the metric $g$ may have
nonpositive scalar curvature $R_g$ (satisfying the above ``tame''
conditions). We consider two cases: $\Lambda >0$ and $\Lambda\leq 0$.
\vspace{2mm}

\noindent 
{\bf Case $\Lambda >0$.}  We have assumed (see above ) that in this
case $R_g>0$ everywhere on $M$, and $R_g\sim \chi^2$ in the
neighborhood of $x_*$. In particular, scalar curvature is positive on
the (compact) belt $B\setminus B_{\sigma }$.  Thus, on $M_{0}\cup
(B\setminus B_{\sigma })$ the scalar curvature is bounded from below
by a positive constant.  Using this we get the estimate
$$ 
C_1 \chi^2 \leq R_g \leq C_2\chi^2
$$ 
everywhere on $M$ (for some positive constants $C_1$, $C_2>0$). We
multiply this inequality by $u^2$, then we integrate over $M$, and
(\ref{c-1}) implies
$$ 
C_1\|u \|_{H^1_2(M)}^2 \leq \<\L_g u, u \> \leq C_2 \|u
\|_{H^1_2(M)}^2
$$ 
for some positive constants $C_1$, $C_2>0$ as required.  
\vspace{2mm}

\noindent 
{\bf Case $\Lambda \leq 0$.}  We use the above decomposition and
(\ref{c-1}) to write 
$$ 
\begin{array}{rcl} \<\L_g u, u \> &= &
\displaystyle \int_{M_{0}\cup B_{\sigma}} \( |\nabla_g u|^2 +
{n-2\over 4(n-1)} R_g u^2 \) d\sigma_g 
\\ 
\\ 
& & \displaystyle 
\ \ \ \ \ \ \ + \int_{(B\setminus B_{\sigma})\cup K} \( |\nabla_g u|^2
+ {n-2\over 4(n-1)} R_g u^2 \) d\sigma_g = A + B.  
\end{array} 
$$ 
We estimate the term $A$: 
$$ 
A\!=\! \int_{M_{0}\cup B_{\sigma}}\! \! \! \!  \( |\nabla_g
u|^2 + {n-2\over 4(n-1)} R_g u^2 \) d\sigma_g \geq C \int_{M_{0}\cup
B_{\sigma}}\! \! \!  \( |\nabla_g u|^2 + \chi^2 u^2 \) d\sigma_g 
$$ 
since $R_g\geq R_0>0$ on $M_0 \cup B_{\sigma}$.  
\vspace{2mm}

For the integral $B$ we have: 
$$ 
\begin{array}{rcl} 
B &= &
\displaystyle \int_{K} \( |\nabla_g u|^2 + {n-2\over 4(n-1)} R_g u^2
\) d\sigma_g 
\\ 
\\ 
& & \displaystyle \ \ \ \ \ \ \ \ \ +
\int_{B\setminus B_{\sigma}} \( |\nabla_g u|^2 + {n-2\over 4(n-1)} R_g
u^2 \) d\sigma_g = {\mathbf I} + {\mathbf J} 
\end{array} 
$$ 
We start with the study of the integrands in the integral ${\mathbf J}$. In the
compact closure $\ov{B\setminus B_{\sigma}}$ we have the norms of the
gradients of a function $u$ with respect to the metrics $g$ and $g_K$:
$$ 
\Vert \nabla
_{g}u\Vert = g^{ij}(x)\xi _{i}\xi_{j},\ \Vert \nabla _{g_{K}}u\Vert =
g^{ij}_{K}(x)\xi _{i}\xi_{j},\ \xi _{i}= \partial _{i}u.  
$$ 
We notice that the ratio
$$ 
\frac{\Vert \nabla _{g}u\Vert _{g}\sqrt{\vert g\vert }}{\Vert
\nabla _{g_{K}}u\Vert _{g_{K}} \sqrt{\vert g_{K} \vert } } =
\frac{g^{ij}_{K}(x)\xi _{i}\xi_{j}\sqrt{\vert g\vert
}}{g^{ij}_{K}(x)\xi _{i}\xi_{j}\sqrt{\vert g_{K}\vert }} 
$$ 
is bounded from below by a positive constant $c$ on the compact set
$\ov{B\setminus B_{\sigma}}$. One can take $c$ to be a minimum of the
ratio taken over the compact subspace 
$$
(x,\xi )\in \{ T(\ov{B\setminus
B_{\sigma}})\ \ \vert \ \ \vert \xi \vert = 1\}
$$ 
of the tangent bundle $T(\ov {B\setminus B_{\sigma}})$.  Thus, on the
compact set $\ov{B\setminus B_{\sigma}}$ we have the bound
\begin{equation}\label{new-1}
\Vert \nabla _{g}u\Vert _{g}\sqrt{\vert g\vert } \geq c\Vert \nabla
_{g_{K}}u\Vert _{g_{K}} \sqrt{\vert g_{K} \vert } 
\end{equation}
for all $u$.  At the same time on this compact set $\ov{B\setminus
B_{\sigma}}$ we have 
 $$ 
\frac{\ell ^{2}\vert R_{g}\vert \sqrt{\vert
g\vert}}{\sqrt{\vert g_{K}\vert }}\leq c_{1} 
$$ 
for some positive constant $c_{1}$ (since the functions $u$ are
continuous and positive on this compact set).  Thus, on this compact belt
we have 
\begin{equation}\label{new-2}
R_{g}\sqrt{\vert g\vert}\geq c_{2}\frac{\Lambda }{\ell
^{2}} \sqrt{\vert g_{K}\vert}, \ \ \ \ \mbox{with some} \ \ c_{2}=
\frac{c_{1}}{\vert \Lambda \vert}>0.  
\end{equation}
Taking $c_{3}= \min (c,c_{2})$ and combining the estimates (\ref{new-1}),
(\ref{new-2}) we get 
$$ 
{\mathbf J}\geq c_{3}\int _{B\setminus B_{\sigma}}\[\Vert \nabla
_{g_{K}}u\Vert ^{2}+\frac{\Lambda }{\ell ^{2}}u^{2}]\sqrt{\vert
g_{K}\vert }\]dx.
$$ 
Thus, 
$$ 
B= {\mathbf I}+{\mathbf J}\geq 
c_{4}\int _{K\cap (B\setminus B_{\sigma})}\[ \Vert \nabla
_{g_{K}}u\Vert ^{2}+\frac{\Lambda }{\ell ^{2}}u^{2}]\sqrt{\vert
g_{K}\vert }\] dx, 
$$ 
here $c_{4}= \min (1,c_{3})>0.$ Thus, we have estimated the integral $B$
from below by an integral over part of the cone, say $\ell \leq
\epsilon _{4}$, containing only \emph{standard conic metric}.
\vspace{2mm}

Thus, we consider the quadratic form 
$$ 
\begin{array}{rcl}
\<\L_{\cone}u,u\>\!\!\! &=&\!\!\! \!\!\!  {\displaystyle
\int _{_{[0,\epsilon ]\!\times\! S^{p}_{r_{p}}\!\times\!
S^{q}_{r_{q}}}}}\!\!\! \! \!\!  \!
\left[\vert \nabla_{\ell } u\vert ^2\! -\!u\!\cdot\! (\Delta
_{\theta}\! + \!\Delta_{\psi }) u \!+\!\frac{n-2}{4(n-1)}\frac{\Lambda }{\ell
^{2}}u^{2}\right]\ell ^{n-1}d\ell d\sigma _{\theta }d\sigma _{\psi }
\end{array}
$$ 
on the space of functions obtained by restriction of functions from
the space $H^{1}_{2}(M)$ to the subset $[0,\epsilon ]\times
S^{p}_{r_{p}}\times S^{q}_{r_{q}}\subset M$ (with the conical standard
metric) and the norm induced from $H^{1}_{2}(M)$.  

Notice that each function from $H^{1}_{2}([0,\epsilon ]\times
S^{p}_{r_{p}}\times S^{q}_{r_{q}})$ can be extended to the function
from $H^{1}_{2}(M)$ (see \cite{Egorov}). Therefore, by restricting
functions from $M$ to the conical part we obtain the whole space
$H^{1}_{2}([0,\epsilon ]\times S^{p}_{r_{p}}\times S^{q}_{r_{q}}).$

Now we decompose the function $u$ into the Fourier series, using the
coordinates $(\ell,\theta,\psi)$ 
$$ 
u= \sum _{i,j}u_{ij}(\ell )\xi
_{i}(\theta )\xi _{j}(\psi ) 
$$ 
and use the notations: 
$$ u_{ij}:=
u_{ij}(\ell), \ \ \ u_{ij,\ell}:= \frac{\p}{\p\ell} u_{ij}(\ell),\ \ \
u_{ij,\ell\ell}:= \frac{\p^2}{\p\ell^2} u_{ij}(\ell).  
$$ 
We have 
$$
\begin{array}{ll} 
&\ \ \ \ \ \ \ -u\Delta_{\psi ,\theta} u 
\\
\\
=&\!\!\!\!\displaystyle 
-\left( \sum _{i,j}u_{ij}\xi _{i}(\theta
)\xi_{j}(\psi )\right)\!\!\cdot \!\! \left( \sum _{i,j}\left[
\frac{2u_{ij}}{r^{2}_{p}\ell ^{2}}\Delta _{\theta }\xi _{i}(\theta )\xi
_{j}(\psi )+
\frac{2u_{ij}\xi _{i}(\theta )}{r^{2}_{q}\ell
^{2}}\Delta _{\psi }\xi _{j}(\psi )\right]\right) 
\\ 
\\ 
=& \!\!\!\!
\displaystyle \left( \sum _{i,j}u_{ij}\xi _{i}(\theta )\xi _{j}(\psi
)\right) \!\!\cdot \!\! \left( \sum _{i,j}\left[- 
\frac{2u_{ij}}{r^{2}_{p}\ell^{2}}\lambda ^{p}_{i}\xi _{i}(\theta )\xi
_{j}(\psi ) -
\frac{2u_{ij}\xi _{i}(\theta )}{r^{2}_{q}\ell
^{2}}\lambda ^{q}_{j}\xi _{j}(\psi )\right]\right) 
\\ 
\\ 
= & \!\!\!\! 
\displaystyle \left(\! \sum_{i,j}u_{ij}\xi _{i}(\theta )\xi _{j}(\psi
)\!\right)\!\!\cdot \!\!  \left(\! \sum _{i,j}\left[
\frac{2}{r^{2}_{p}\ell ^{2}}\lambda
^{p}_{i}\!+\!\frac{2}{r^{2}_{q}\ell ^{2}}\lambda ^{q}_{j}\right]
u_{ij}(\ell )\xi _{i}(\theta )\xi _{j}(\psi )\!\right).  
\end{array}
$$ 
Now we add the term 
$$
\vert \nabla _{r}u\vert
^2+\frac{n-2}{4(n-1)}\frac{\Lambda }{\ell^{2}}u^{2},
$$ 
to the last expression, where $u$ is decomposed into the same Fourier
series by $\psi, \theta $. In particular, we have
$$
\vert \nabla_{\ell} u\vert ^2 = \vert \sum
_{i,j}[\nabla _{\ell }u_{ij}(\ell )]\xi _{i}(\theta )\xi _{j}(\psi )
\vert ^2 .  
$$ 
Then we integrate the resulting expression over the
product of spheres $S^{p}\times S^{q}$. We obtain 
$$ 
\<\L_{\cone}u,u\>= \int _{[0,\epsilon ]}\ell^{n-1} d\ell \sum
_{ij}\left[ u_{ij,\ell}^{2}+ K_{ij}u_{ij}^{2} \right] = \sum
_{ij}\int^{\epsilon }_0\ell ^{n-1}d\ell [u^{2}_{ij,\ell}+K_{ij}u^{2}]
, 
$$ 
where 
$$ 
K_{ij}= \frac{2 \lambda ^{p}_{i}}{r^{2}_{p}}+
\frac{2\lambda ^{q}_{j}}{r^{2}_{q}}+\frac{n-2}{4(n-1)}\Lambda .  
$$
The total quadratic form $\<\L u,u\>$ is estimated from below by the
integral over $M_{0}$ in $\<\L u,u\>$ plus some positive constant
times the form $\<\L_{\cone }u,u\>$: 
$$ 
\<\L u,u\> \geq  c\<\L_{\cone}u,u\>+
\int _{M_{0}}[\vert \nabla u\vert^2 +R_{g}u^{2} ]d\sigma _{g}.
$$ 
Since $R_{g}>0$ on $M_{0}$, second term is always positive. This is
not true for the first term. Take for example $u= u_{0}= const$. This
function belongs to the space $H^{1}_{2}(M)$ and $\<\L_{\cone
}u_{0},u_{0}\>= K_{00}u_{0}^{2}\frac{\epsilon ^{n}}{n}$, thus is
negative for $\Lambda <0 $, since $\lambda^{p}_{0}= \lambda ^{q}_{0}=
0.$ This leads, therefore, to the following necessary condition for
the form $\<\L u,u\>$ to be positive (take $u= 1$):
$$ 
\int _{M}R_{g}d\sigma _{g}>0.  
$$ 
Now we find a lower bound for the quadratic form values of 
$$ \int_0^{\epsilon }\left[\ell^{n-1} v^{2}_{,\ell }
+ K \ell^{n-3}v^{2}\right] d\ell 
$$ 
on the space $H^{1}_{2}(0,\epsilon)$.  Later we will specialize the
results to the cases $v= u_{ij},K= K_{ij}$.
\vspace{2mm}

Denote by $\bar v$ the function on $(0,1)$ obtained by the scaling
$\ell = \epsilon _{3}t$: $v(\ell )= v(\epsilon _{3}t)= {\bar v}(t).$
We have 
$$ 
\int_{0}^{\epsilon }[v^{2}_{,\ell } + Kv^{2}\ell ^{-2}]\ell
^{n-1}d\ell = \epsilon ^{n-2}\int_{0}^{1}[{\bar v}_{,t}^{2}+K{\bar
v}^{2}t^{-2}]t^{n-1}dt.
$$ 
Change of variables $t= s^{a}$; $dt= as^{a-1}ds$; $s= t^{1\over a}$ 
gives 
$$
f_{,s}= f_{,t}t_{,s}= f_{,t}as^{a-1}; \ \ \ f_{,t}=
\frac{1}{a}s^{1-a}f_{,s}.
$$ 
This  leads to
$$  
\begin{array}{rcl}  
\!\!\displaystyle \int _{0}^{\epsilon}[v^{2}_{,\ell }+Kv^{2}\ell
^{-2}]\ell ^{n-1}d\ell \!\!\!\!&= &\!\!\!\! \displaystyle 
\epsilon^{n-2}\!\!\int\!\!
s^{a(n-1)}as^{a-1}\left[\frac{1}{a^2}s^{2(1-a)}{\bar
v}^{2}_{,s}+Ks^{-2a}{\bar v}^{2}\right] ds 
\\ 
\\ 
&\!\!= \!\!& \!\!\!\!\displaystyle
\epsilon ^{n-2}a\int s^{a(n-1)+a-1+2-2a}\left[\frac{1}{a^2}{\bar
v}^{2}_{,s}+Ks^{-2}{\bar v}^{2}\right] ds  
\\
\\
&\!\!= \!\!& \!\!\!\!\displaystyle
-(n-2)\epsilon ^{n-2}\int ^{1}_{+\infty }\[{\bar
v}^{2}_{,s}+a^{2}Ks^{-2}{\bar v}^{2}\]ds 
\\
\\
&\!\!= \!\!& \!\!\!\!\displaystyle
(n-2)\epsilon ^{n-2}\int
_{1}^{+\infty }\[{\bar v}^{2}_{,s}+\frac{K}{(n-2)^{2}}s^{-2}{\bar
v}^{2}\]ds.
\end{array} 
$$ 
Here we took $a= -\frac{1}{n-2}$.  Thus it is enough to give a lower
bound for
$$ 
\int _{1}^{+\infty } \[{\bar v}^{2}_{,s}+K_{1}{\bar
v}^{2}(s)\] ds, \ \ \ K_{1}= \frac{K}{(n-2)^{2}}, 
$$ 
on the space of $H^{1}_{2}(1,+\infty )$ (with the norm
defined by the same quadratic expression with $K_{1}= 1$).  
\vspace{2mm}

Functions from Sobolev space $H^{1}_{2}(1,+\infty )$ are continuous at
$s= 1$ and have the well-defined limit value $f(1)= \lim_{s\to
1}f(s)$.  This value is the continuous linear functional on the space
$H^{1}_{2}(1,+\infty )$. The kernel of this functional is the space
$H^{1}_{2,0}(1,\infty)$ (which is the closure in
$H^{1}_{2,0}(1,\infty)$ of the subspace $C_{0}^{\infty }(1,\infty )$
of smooth functions with compact support \cite{Egorov}). A function
$f$ from $H^{1}_{2,0}$ has a canonical extension $\hat f$ by zero to
the function in $H^{1}_{2}(0,+\infty )$ with the norm defined by the
same quadratic form
$$ 
\|{\hat f}\|^{2}= \int _{0}^{+\infty } \[{\hat
f}^{2}_{,s}+{\hat f}^{2}(s)\]ds.  
$$ 
Notice that the extension $\hat f$ has the same norm as the function
$f$.  For such a function, obtained by the extension to $(1, \infty )$
of a function from $H^{1}_{2,0}(1, \infty )$ and thus, being zero in a
neighborhood of zero, we can use the simplest Hardy inequality (see
\cite{HL}) to estimate
$$ 
\int_{0}^{+\infty }{\hat f}^{2}(s)s^{-2}ds\leqq 4\int _{0}^{+\infty }{\hat
f}_{,s}^{2}(s)ds.  
$$ 
Because of the construction of the extension we get
a similar inequality with $f$ instead of $\hat f$ and the lower
limit $1$ replacing zero.  
\vspace{2mm}

Using this estimate we get for $f\in H^{1}_{2,0}(1,\infty)$ and
negative $K_1$: 
$$
\int _{1}^{+\infty } \[f^{2}_{,s}+K_{1}f^{2}(s)\]ds
\geq (1+4K_{1})\int _{1}^{+\infty } f^{2}_{,s}ds .  
$$ 
Thus, if $1+4K_{1} >0$, the quadratic form 
$$ 
K(f,f)= \int
_{1}^{+\infty }\[f^{2}_{,s}+K_{1}f^{2}(s)\] ds 
$$ 
is positive definite on $H^{1}_{2,0}(1,\infty)$, where the norm
induced from $H^{1}_{2}(1,\infty)$, is equivalent to the norm $\int
_{1}^{+\infty } f^{2}_{,s}ds $ (this follows from the Hardy
inequality). Applying this to the case where $K= K_{ij}$ we see that
it is sufficient to check this condition for the case $i= j= 0$ (since
$\lambda ^{p}_{i},\lambda ^{q}_{j} >0$).  For $i= j= 0$,
$$ 
K_{1}= \frac{K}{(n-2)^{2}}=
\frac{\Lambda
}{4(n-1)(n-2)}.  
$$ 
Therefore
$$ 
1+4K_{1}= 1+4\frac{\Lambda }{4(n-1)(n-2)}= 1+\frac{\Lambda
}{(n-1)(n-2)}= \mu ^{2}>0, 
$$ 
see Appendix.  Therefore, on the subspace $H^{1}_{2,0}(1,\infty)$ the
quadratic form $K(f,f)$ is positive definite.  
\vspace{2mm}

To determine what happens at the complement to this
subspace and to estimate our quadratic form from below we return to
the quadratic form on the space $H^{1}_{2}(0,\epsilon)$ with the
(square) of the norm 
$$ 
\int _{0}^{\epsilon}
\ell^{n-1}(v_{,\ell}^{2}+\ell ^{-2}v^{2})d\ell = \Vert v_{,\ell}\Vert
^{2}_{L_{2}(0,\epsilon;\ell ^{n-1}d\ell) }+\Vert \ell ^{-1}v\Vert
^{2}_{L_{2}(0,\epsilon;\ell ^{n-1}d\ell)} .
$$   
Let $K\leqq 0$. To find the lowest eignevalue of the quadratic form 
$$
\int _{0}^{\epsilon} \ell^{n-1}(v_{,\ell}^{2}+K\ell ^{-2}v^{2})d\ell
$$ 
we calculate the minimum of the following fraction 
$$ 
\min_{u\ne
0}\frac{\Vert u_{,s}\Vert _{L_{2}}+K \Vert s^{-1}u\Vert
_{L_{2}}}{\Vert u_{,s}\Vert _{L_{2}}+\Vert s^{-1}u\Vert _{L_{2}}} 
$$
over $u\in H^{1}_{2}(0,\epsilon)$.  We write this relation as 
$$
\frac{f(u)+K}{f(u)+1}, \ \ \ \mbox{with} \ \ \ f(u)= \frac{\Vert
v_{,l}\Vert ^{2}_{L_{2}(0,\epsilon;\ell ^{n-1}d\ell }}{\Vert \ell
^{-1}v\Vert ^{2}_{L_{2}(0,\epsilon;\ell ^{n-1}d\ell }}.  
$$ 
Now we notice that the function 
$$ 
s\rightarrow \frac{s+K}{s+1} 
$$ 
is increasing for $K \leqq 0$ and for $s\geqq 0$ takes its minimal
value (equal to $K$) at $s= 0$. On the other hand $f(u)$ is well
defined for all $u\in H^{1}_{2}(0,\epsilon)$, $u\neq 0$, (since its
denominator cannot be zero for $u\ne 0$). The function $f(u)$ is
nonnegative and is equal to zero only if $u_{,\ell }= 0$, i.e. for
constant functions $u= const$.
\vspace{2mm}

Therefore, the quadratic form 
$$ 
\int _{0}^{\epsilon} \ell^{n-1}(v_{,\ell}^{2}+K\ell ^{-2}v^{2})d\ell
$$  
has $K$ as its minimal eigenvalue and the constant $u_{0}=
\sqrt{\frac{n-2}{\epsilon ^{n-2}}}$ as its eigenvector of unit norm,
corresponding to this eigenvalue.  
\vspace{2mm}

As a result, the condition $\Lambda >0$, or, what is the same, $\mu
>1$ is sufficient for the total quadratic form $L$, and the operator
$\L_{g}$ to be positive definite.  Recall that we have 
$$
\begin{array}{rcl} \<\L u,u\>& \geq & c\<\L_{\cone }u,u\>+\int
_{M_{0}}[\vert \nabla u\vert +R_{g}u^{2} ]d\sigma _{g} 
\\ 
\\ 
&= &\displaystyle \sum _{ij}\int ^{\epsilon }\ell ^{n-1}
[u^{2}_{ij}+K_{ij}u^{2}]+\int _{M_{0}}[\vert \nabla u\vert +R_{g}u^{2}
]d\ell d\sigma _{g}.  
\end{array} 
$$ 
The integral over $M_{0}$ is always nonnegative.  In the sum, terms
for which $K_{ij}\geq 0$ are also nonnegative.  Since 
$$ 
K_{ij}=
\frac{2 \lambda ^{p}_{i}}{r^{2}_{p}}+\frac{2\lambda
^{q}_{j}}{r^{2}_{q}}+\frac{n-2}{4(n-1)}\Lambda = \frac{2 \lambda
^{p}_{i}}{r^{2}_{p}}+\frac{2\lambda
^{q}_{j}}{r^{2}_{q}}+\frac{(n-2)^{2}(\mu^{2}-1)}{4}, 
$$ 
for $p,r_{p},q,r_{q}$ fixed, all the terms in the sum over $i,j$ are
nonnegative except a finite number of them.  
\vspace{2mm}

More than this, each term in the sum (quadratic form) for which
$K_{ij}\leqq 0$ has one and only one nonpositive eigenvalue with
constant eigenfunction.  For the form $\L_{\cone}$ on
$H^{1}_{2}((0,\epsilon)\times S^{p}\times S^{q})$ this corresponds to
the function(s) $u_{ij}= c_{ij}\xi ^{p}_{i}(\theta )\xi^{q}_{j}(\psi)$
with some constants $c_{ij}$. The condition $K_{ij}\leqq 0$ can be
rewritten as follows
$$
\frac{ \lambda^{p}_{i}}{r^{2}_{p}}+\frac{\lambda
^{q}_{j}}{r^{2}_{q}}\leqq \frac{(n-2)}{(n-1)}\vert \Lambda \vert
$$ 
and the number of negative modes of the form $\L_{\cone}$ can be found
from here. Notice, in particular, that if $\Lambda <0$, then for $i,j=
0$ that condition is always satisfied; thus, there is at least one
negative mode for $\L_{\cone}$.  That does not prevent, though, this
negative input in the whole form $\<\L u,u\>$ from being compensated
by the input of the $M_{0}$-part.
\vspace{2mm}

We also notice that as it follows from the proof, in the case where
$\Lambda \leqq 0$, there is the lower bound for the form $\<\L u,u\>$,
i. e.
$$
\<\L u,u\>\geqq D\Vert u\Vert _{H^{1}_{2}(M)}.
$$ 
Here $D$ is a finite constant, and 
$$
|D| \leq \frac{4(n-2)}{n-2}\vert \Lambda \vert c_{4},
$$
where $c_{4}$ is the constant above.
This proves Theorem
\ref{Th1new}.
\end{pf-of} 
{\bf \ref{s4}.4. Necessary condition for positivity of $\<\L u,u\>$.}
We examine the case when $\Lambda\leq 0$ but $K_{ij}>0$ for all $i>0$
and $j>0$. Since
$$
\Lambda = - (n-1)(n-2) + \frac{2p(p-1)}{r^2_p} + \frac{2q(q-1)}{r^2_q},
$$
then the condition $\Lambda\leq 0$ is equivalent to 
\begin{equation}\label{neg_1}
\frac{(n-1)(n-2)}{2} \geq \frac{p(p-1)}{r^2_p} + \frac{q(q-1)}{r^2_q} .
\end{equation}
On the other hand the conditions $K_{10}>0$ and $K_{01}>0$ imply that
$K_{ij}>0$ for all $i>0$ and $j>0$. We recall that zero has
multiplicity $1$ for $\Delta_{\theta}$ and $\Delta_{\psi}$, and the next
eigenvalue is $p$ for $\Delta_{\theta}$ and $q$ for $\Delta_{\psi}$
respectively. This gives that $K_{10}>0$ is equivalent to
$$
\begin{array}{l}
\displaystyle
\frac{2p}{r_p^2} + \frac{(n-2)^2\Lambda}{4(n-1)(n-2)} > 0 \ \ \mbox{or}
\ \ 
\Lambda > -\frac{8(n-1)}{n-2}\cdot \frac{p}{r_p^2}.
\end{array}
$$
We use the above formula for $\Lambda$ to get that
\begin{equation}\label{neg_2}
K_{10}>0 \ \ \Longleftrightarrow \  \ 
\frac{4(n-1)}{n-2} \cdot \frac{p}{r_p^2}
+ \frac{p(p-1)}{r^2_p} + \frac{q(q-1)}{r^2_q} > \frac{(n-1)(n-2)}{2} 
\end{equation}
Similarly, 
\begin{equation}\label{neg_3}
K_{01}>0 \ \ \Longleftrightarrow \  \ 
\frac{4(n-1)}{n-2} \cdot \frac{q}{r_q^2}
+ \frac{p(p-1)}{r^2_p} + \frac{q(q-1)}{r^2_q} > \frac{(n-1)(n-2)}{2} 
\end{equation}
We have proved the following.
\begin{Proposition}\label{neg_4}
Assume $\Lambda\leq 0$ {\rm (}which is equivalent to {\rm
(\ref{neg_1})}{\rm )} and the conditions {\rm (\ref{neg_2})}, {\rm
(\ref{neg_3})} are satisfied. Then the form $\<\L u,u\>$ is positive
if and only if
$$
\int_M R_g d\sigma_g > 0.
$$
\end{Proposition}
\section{Weak Yamabe Theorem}\label{s5} 
{\bf \ref{s5}.1. Yamabe functional.}  Now we define the Yamabe
functional on $M$ and study its properties.  Let $\alpha\in
[1,\alpha^*]$, where $\alpha^*=  {n+2\over n-2}$. For each
$\alpha$ we consider the functional
$$
\begin{array}{c}
\displaystyle I_{\alpha}(\phi)=  {E(\phi)\over \(\int_M
|\phi|^{\alpha+1} dV_g\)^{2\over \alpha+1}}, \ \ \ \ \phi\in
H^1(M,dV_g),\ \ \phi \ne 0.
\\
\\
\displaystyle E(\phi) =  \int\( |\nabla \phi|^2 + {n-2\over 4(n-1)} R_g
\phi^2\) dV_g .
\end{array}
$$
This is the Yamabe functional if $\alpha= \alpha^*$. One can prove the
following fact by using Theorem \ref{sob-6}.
\begin{Proposition}\label{L2}
The functional $I_{\alpha} : H_2^1(M,g)\longrightarrow
\R $ is defined and continuous on the space $H_2^1(M,g)$ for all
$\alpha\in [1,\alpha^*)$.
\end{Proposition}
We denote by $C_{\mu}$ the {\sl norm of the embedding} $H_2^1(M) \subset
L_{\mu}(M)$, i.e.
$$
C_{\mu}= 
\inf_{^{\begin{array}{c}\phi\neq 0\\\phi\in H_2^1(M)\end{array}}}
\frac{\|\phi\|_{L_{\mu}}}{\|\phi\|_{H^1_2}},
$$
where it is assumed that $\mu\leq 1+\alpha^*$.
\begin{Proposition}\label{Prop_new1}
Let $\Lambda >0$. Then the functional $I_{\alpha}$ is bounded from
below, i.e. for all $\phi\in H_2^1(M)$, $\phi\neq 0$  we have
$$
I_{\alpha}(\phi) \geq \frac{C_1}{C_{\alpha+1}^2}
$$
with the constant $C_1>0$ given in Theorem \ref{Th1new}, and
$C_{\alpha+1}$ as above.
\end{Proposition}
\begin{Proof}
Theorem \ref{Th1new} gives that $\<\L_gu,u\> \geq C_1\|u\|^2_{H^1_2}$.
Thus
$$
I_{\alpha}(\phi) = 
\frac{\<\L_g\phi,\phi\>}{\|\phi\|_{L_{\alpha+1}}^2} \geq
\frac{C_{1}\|\phi\|_{H^1_2}^2}{\|\phi\|_{L_{\alpha+1}}^2}
$$
for any $\phi\in C^{\infty}_*(M)$. The embedding Theorem
\ref{sob-6} gives the continuous embedding $H^1_2(M)\subset
L_{\alpha +1}(M)$ with $\alpha\leq \alpha^*$, and
$\|\phi\|_{L_{\alpha+1}}\leq C_{\alpha+1}\|\phi\|_{H^1_2}$. Thus
$$
\frac{\|\phi\|_{H^1_2}^2}{\|\phi\|_{L_{\alpha+1}}^2} \geq 
\frac{1}{C_{\alpha+1}}
\ \ \ \mbox{giving} \ \ \
I_{\alpha}(\phi) \geq \frac{C_1}{C_{\alpha+1}^2}.
$$
Since the space $C^{\infty}_*(M)$ is dense in $H^1_2(M)$, it gives the
result.
\end{Proof}
Here is the ``easy'' version of the Yamabe theorem. To prove it we
follow the course of the corresponding result for closed
manifolds, see \cite[Theorem 5.5]{Aubin3}.
\begin{Theorem}\label{sob-11}
Let $M$ be a compact closed manifold with a metric $g$ and a
conical singularity as above. Let $\Lambda >0$. For any $\alpha\in
[1,\alpha^*)$ there exists a function $u_{\alpha}\geq 0$
minimizing the functional $I_{\alpha}$, so that
$$
\begin{array}{c}
\displaystyle
\int_M u_{\alpha}^{\alpha+1} dV_g = 1
\\
\\
\displaystyle \mbox{giving} \ \ \ \ I_{\alpha}(u_{\alpha}) =  \min
\{ I_{\alpha}(\phi) \ | \ \phi\in H^1_{2}(M),\ \phi \ne 0 \}.
\end{array}
$$
\end{Theorem}
Denote this value $Q_{\alpha}= 
I_{\alpha}(M,[g])= I_{\alpha}(u_{\alpha})$. The function
$u_{\alpha}$ is a weak (in $ H^{1}_{2}(M)$) solution of the
equation:
$$
 -\Delta u_{\alpha}+ {n-2\over 4(n-1)} R_{g}u_{\alpha} = Q_{\alpha}
u_{\alpha}^{\alpha}.
$$
{\bf \ref{s5}.1. Proof of Theorem \ref{sob-11}.}  {\bf (a)} First we
prove that the functional $I_{\alpha}$ is bounded, and thus
$Q_{\alpha}$ is finite. Let $q= 1+\alpha$, $2\leq q < 1+\alpha^*$.
The conformal Laplacian $\L_g$ is positive by Theorem \ref{Th1new}
(Recall that $\Lambda >0$). Thus $I_{\alpha}(\phi)\geq 0$ for any
$\phi\in H_2^1(M,g)$, therefore $Q_{\alpha}\geq 0$ (moreover,
$Q_{\alpha }>\frac{C_{1}}{C^{2}_{\alpha +1}}$, see Proposition
\ref{Prop_new1}). On the other hand,
$$
\begin{array}{rcl}
Q_{\alpha} \leq I_{\alpha}(1) 
& = & 
\displaystyle
{n-2\over 4(n-1)}\int_M R_g d\sigma_g
\\
\\
& = & 
\displaystyle
{n-2\over 4(n-1)}\int_{M\setminus K} R_g d\sigma_g  +
{n-2\over 4(n-1)}\int_K R_g d\sigma_g .
\end{array}
$$
The first integral on the right is bounded since $R_g$ is
continuous in $M\setminus \{x^*\}$. The second integral is
bounded since $n>2$ and
$$
R_r\sim {C\over \ell^2} \ \ \mbox{as} \ \ \ell \rightarrow 0, \ \
\mbox{and} \ \ \int_{0}^{\epsilon_0}\ell^{n-3}d\ell <\infty.
$$
{\bf (a$^{\prime}$)} Now we have that $\Vol_g(M)<\infty$, and
$1\in L_s(M)$ for any $s\geq 2$. Thus we use the inequality
$$
\int_M fg d\sigma_g \leq \(\int_M f^s d\sigma_g \)^{1/s} \(\int_M
g^{s'} d\sigma_g \)^{1/s^{\prime}}
$$
(which holds for positive functions $f\in L_s(M)$, $g\in
L_s^{\prime}(M)$, with ${1\over s }+ {1\over s^{\prime}} = 1$,
$s,s^{\prime}>0$).  We apply this for $f= \phi^q$, $g= 1$,
$s= \frac{q}{2}$, $s^{\prime}= \frac{q}{q-2}>0$.
Thus we get
$$
\begin{array}{l}
\displaystyle
\int_M \phi^2\cdot 1 d\sigma_g \leq
\(\int_M (\phi^2)^{\frac{q}{2}}d\sigma_g \)^{\frac{2}{q}} \cdot
\(\int_M 1^{\frac{q}{q-2}}d\sigma_g \)^{\frac{q-2}{q}}, \ \ \ \mbox{or}
\\
\\
\displaystyle
\|\phi\|_{L_2(M)}^2\leq \|\phi\|_{L_q(M)}^2\cdot 
\(\Vol_g(M)\)^{\frac{q-2}{q}},
\ \ \ \mbox{so}
\\
\\
\displaystyle
\|\phi\|_{L_2(M)}\leq \|\phi\|_{L_q(M)}\cdot 
\(\Vol_g(M)\)^{\frac{q-2}{2q}}.
\end{array}
$$
{\bf (b)} Now let $\{\phi_i\}$ be a minimizing sequence such that
$$
\int_M\phi_i^qd\sigma_g =  1, \ \ \phi_i\in H_2^1(M,g), \ \ \mbox{and} 
\ \
\lim_{i\rightarrow\infty} I_{\alpha}(\phi_i)=  Q_{\alpha}.
$$
First we prove that the set $\{\phi_i\}$ is bounded in
$H_2^1(M,g)$. We have
\begin{equation}\label{sob-12}
\begin{array}{rcl}
\displaystyle
\|\phi_i\|_{H_2^1(M)} &= & \|\nabla\phi_i\|_{L_2(M)}^2 +
\|\chi\phi_i\|^2_{L_2(M)}
\\
\\
 &= & \displaystyle
I_{\alpha}(\phi_i)-{n-2\over
4(n-1)}\int_{M}R_g\phi_i^2d\sigma_g +\int_{M}\chi^2\phi_i^2d\sigma_g .
\end{array}
\end{equation}
Since $\{\phi_i\}$ is a minimizing sequence, we can assume that
$I_{\alpha}(\phi_i)\leq Q_{\alpha}+1$.
\vspace{2mm}

\noindent Now we consider the case when $R_g>0$ everywhere (i.e.
$\Lambda>0$). Then we have that $R_g= {\Lambda\over \ell^2}$ on
the cone $K$, thus $\chi^2(x) < C R_g(x)$ for some positive
constant $C$ and any $x\in M$. The sum of the first two terms in
(\ref{sob-12}) coincides with $\|\nabla\phi_i\|_{L_2(M)}^2$, so
it is positive. Therefore
$$
\|\phi_i\|_{H_2^1(M)} \leq A \cdot I_{\alpha}(\phi_i)-A \cdot {n-2\over
4(n-1)}\int_{M}R_g\phi_i^2d\sigma_g + C \int_{M}R_g\phi_i^2d\sigma_g
$$
for any $A\geq 1$. We choose $A$ large enough, so that
$$
C - A {n-2\over 4(n-1)} < 0,
$$
to get the estimate
$$
\|\phi_i\|_{H_2^1(M)} \leq A\cdot I_{\alpha}(\phi_i) \leq
A\cdot (Q_{\alpha}+1).
$$
Notice that on $M\setminus K$ both integrals are estimated by the norm
$\|\phi\|_{L_2(M)}$.
\vspace{2mm}

\noindent
{\bf (c)} We follow the proof of {\bf (c)} in \cite[Theorem
5.5]{Aubin3} to find a subsequence $\{\phi_j\}$ of $\{\phi_i\}$ and a
nonnegative function $u_{\alpha}\in H^1_2(M)$ such that
\begin{enumerate}
\item[$(\alpha)$] $\phi_j \to u_{\alpha}$ in $L_{\alpha+1}(M)$;
\item[$(\beta)$]  $\phi_j \to u_{\alpha}$ weakly in $H^1_2(M)$;
\item[$(\gamma)$] $\phi_j \to u_{\alpha}$ almost everywhere.
\end{enumerate}
One may satisfy $(\gamma)$ since
$L_{\alpha+1}(M)\subset L_2(M)$ continuously and any sequence
converging in $L_2(M)$ has a subsequence that converges almost
everywhere. To satisfy $(\beta)$ one uses that the embeddings
$H_1^2(M)\subset L_{\alpha+1}(M)\subset L_2(M)\subset H_2^1(M)$ are
continuous. Then $u_{\alpha}\in H_1^1(M)$ because of weak compactness
in reflexible Banach spaces, see \cite[Chapter V, Section 2]{Yosida}.
Finally \cite[Chapter V, Section 1]{Yosida} gives
$$
\|u_{\alpha}\|_{H^1_1}\leq \lim_{\ \ \ j\to\infty}\!\!\!\!{\mathrm =
i}{\mathrm n}{\mathrm f} \|\phi_j\|_{H^1_1}.
$$
This proves {\bf (c)}.
\vspace{2mm}

\noindent
{\bf (d)} Here we prove that $u_{\alpha}$ is a weak solution of the
Yamabe equation. It means that for all $\phi\in H_2^1(M)$
$$
\int_M (\nabla u_{\alpha})\cdot(\nabla \phi)d\sigma_g +
\frac{n-2}{4(n-1)}\int_M R_g u_{\alpha}\phi d\sigma_g =  
Q_{\alpha}\int_{M}
u_{\alpha}^{\alpha} \phi d\sigma_g.
$$
The proof is literally the same as in \cite[Theorem 5.5]{Aubin3}.
The only difference is the use of the space $C_*^{\infty}(M)$
instead of $C_0^{\infty}(M)$. This ends the proof of theorem
5.3 for $\lambda >0$.
 \hfill $\Box$
\begin{Corollary}
The solution $\phi_{\alpha}\neq 0$ on $M_*$.
\end{Corollary}
{\bf Remark.} 
{\rm We will see later, in Secion \ref{s7}, that there are cases
where the minimizer (Yamabe solution) belongs to the space
$H^{1}_{2}(M)$ but not to $H^{2}_{2}(M)$. This is very different from
the case of compact manifolds (cf. \cite{Aubin1}).}
\vspace{2mm}

\noindent
{\bf Remark.} 
{\rm The case $\Lambda\leq 0$ is also very interesting. One
can prove the result similar to Theorem \ref{sob-11} under the same
restrictions as in Proposition \ref{neg_4}.}
\section{Asymptotic of solutions: the linear case}\label{s6}
In this section we study the asymptotic behavior of solutions of the
linear equation $\L_{g}u=  Q_{1}u$ near the point $x_*$.  We use the
polar coordinates $(\ell, \theta,\phi )$ on the conical part $K$. Then
the equation $\L_{g}u=  Q_{1}u$ has the form
\begin{equation}\label{cone1-10}
\begin{array}{l}
\( {\p\over \p\ell^2} + {n-1\over \ell}{\p\over \p\ell}+ {2\over
{r_p^{2}\ell^2}}\Delta_{\theta}+{2\over
{r_q^{2}\ell^2}}\Delta_{\psi}\)u+
\[ -{n-2\over 4(n-1)} {\Lambda\over \ell^2} + Q_1\]u =   0.
\end{array}
\end{equation}
First, we recall some basic information on the Laplacian operator on
spheres (\cite{Gilkey}). Let $\{\lambda_j^p,\chi_j^p\}$ be the
spectrum of the Laplacian $\Delta_{S^p}= -\Delta_{\theta}$, and,
respectively, $\{\lambda_i^q,\chi_i^q\}$ of $\Delta_{S^q}=
-\Delta_{\psi}$.  It is well-known (\cite{Gilkey}) that
$$
\begin{array}{l}
\lambda_0^p =   0, \ \ \lambda_1^p =   \cdots \lambda_{p+1}^p =  p,
\ \ \lambda_{p+2}^{p}=  2(p+1),\ldots ,
\\
\\
\lambda_0^q =   0, \ \ \lambda_1^q =   \cdots \lambda_{q+1}^q =  q,
\ \ \lambda_{q+2}^{q}=  2(q+1),\ldots .
\end{array}
$$
Any $L^2$-function on $S^p$ (correspondingly on $S^q$) decomposes into
Fourier series with respect to the orthonormal basis $ (\chi_i^p)$
(correspondingly $(\chi_j^q )$). On the cone $K$ we have
\begin{equation}\label{cone1-9}
u(\ell,\theta,\psi) =   \sum_{ij} u_{ij}(\ell)
\chi_i^p(\theta)\chi_j^q(\psi).
\end{equation}
We decompose $u$ as in (\ref{cone1-9}) to obtain the following system
of equations for the coefficient functions $u_{ij}(\ell)$ on the 
half-line $\ell \geq 0$:
\begin{equation}\label{cone1-11}
\begin{array}{c}
{\p^2 u_{ij}\over \p\ell^2} + {n-1\over \ell} {\p u_{ij}\over \p \ell}
+ \[ Q_1 - \( {n-2\over 4(n-1)}\Lambda +\frac{2}{r_p^{2}}\lambda_i^p +
\frac{2}{r_q^{2}}\lambda_j^q\){1\over\ell^2}\] u_{ij} =  0, \ \ \ 
\mbox{or}
\\
\\
{\p^2 u_{ij}\over \p\ell^2} + {n-1\over \ell} {\p
u_{ij}\over \p \ell} + \[ Q_1 - {K_{ij}\over \ell^2}\] u_{ij}=  0.
\\
\\
\mbox{with} \ \ K_{ij} =   {n-2\over 4(n-1)}\Lambda +
\frac{2}{r_p^{2}}\lambda_i^p + \frac{2}{r_q^{2}}\lambda_j^q.
\end{array}\!\!\!\!\!\!\!\!\!
\end{equation}
The equations (\ref{cone1-11}) are known as {\sl degenerate
hypergeometric} or {\sl Whitteker equations} (see \cite[Vol. 1,
Chapter 6]{Bateman}).  Such an equation can be reduced, via an
appropriate substitution, to the Bessel equations with the pure
imaginary parameter $\nu $.  Their solutions can be explicitly
written in terms of the corresponding Bessel functions.
\vspace{2mm}

\noindent
Here we are interested in asymptotic behavior of solutions as $\ell
\rightarrow 0$. Thus, we are looking for solutions in the form of
power series
$$
u_{ij} =   \ell^{\nu_{ij}}\sum_{k=  0}^{\infty} a_k \ell^k .
$$
We have the first and the second derivatives:
$$
\begin{array}{rcl}
\!\!u_{ij}^{\prime}\! &\!\! =  \!\! &\! \displaystyle
\nu_{ij}\ell^{\nu_{ij}-1}\sum_{k=  0}^{\infty} a_k \ell^k
+ \ell^{\nu_{ij}}\sum_{k=  0}^{\infty} (k+1) a_{k+1} \ell^k,
\\
\\
\!\!u_{ij}^{\prime\prime}\!\! & \!=  \! \!& \!\displaystyle
\nu_{ij}(\nu_{ij}-1)\ell^{\nu_{ij}-2}\!\sum_{k=  0}^{\infty} a_k \ell^k
\!+ \!2 \nu_{ij} \ell^{\nu_{ij}-1}\!\sum_{k=  0}^{\infty} (k+1) a_{k+1} 
\\
\\
 & \!=  \! \!& \!\displaystyle
\ell^k
\!+ \!\ell^{\nu_{ij}}\!\sum_{k=  0}^{\infty})\! (k+2)(k+1)a_{k+2} \ell^k .
\end{array}
$$
We collect the coefficients for the different powers of $\ell$ in
the equation (\ref{cone1-11}):
$$
\begin{array}{rcl}
\ell^{\nu_{ij}-2} & : & \displaystyle \nu_{ij}(\nu_{ij}-1)a_0 +
(n-1)\nu_{ij} a_0 - K_{ij} a_0 =   0,
\\
\\
\ell^{\nu_{ij}-1} & : & \displaystyle \nu_{ij}(\nu_{ij}-1)a_1 + 2
\nu_{ij} a_1 + (n-1)(\nu_{ij} + 1) a_1 -K_{ij} a_1=  0.
\end{array}
$$
For the (general) coefficient of $\ell ^{\nu_{ij}+m}$ we obtain
the following equation:
$$
\begin{array}{rcl}
\ell^{\nu_{ij}+m} & : & \displaystyle
\nu_{ij}(\nu_{ij}-1)a_{m+2}  + 2\nu_{ij}(m+2) a_{m+2} +
(m+2)(m+1)a_{m+2}
\\
\\
&  & +
 (n-1)(\nu_{ij}+m+2)a_{m+2}
+Q_1 a_m - K_{ij}a_{m+2}=  0.
\end{array}
$$
Thus we get the recursive equation for the coefficients $a_{m}$ which
we denote by $(Y_{m+2})$:
\vspace{2mm}

\noindent \hspace*{5mm} $[(\nu_{ij}+m+2)^{2} +
(n-2)(\nu_{ij}+m+2)-K_{ij}] a_{m+2} =  - Q_1 a_m$ \hfill
$(Y_{m+2})$ \vspace{2mm}

\noindent
Denote by $K_{\nu_{ij}+m}$ the left side of the
previous equation. We consider the equation $(Y_0)$:
\vspace{2mm}

\noindent \hspace*{30mm} $(\nu_{ij}^2 + (n-2)\nu_{ij} - K_{ij})a_0
=   0$.\hfill $(Y_{0})$ \vspace{2mm}

\noindent
Here we have either $a_0=  0$ or
$$
\nu_{ij}^{(e,\pm)} =   -{n-2\over 2}\pm \sqrt{\({n-2\over 2}\)^2 +
K_{ij}}.
$$
The second equation $(Y_1)$ is as follows:
\vspace{2mm}

\noindent \hspace*{30mm}  $(\nu_{ij}^2 + n\nu_{ij} +(n-1 - K
_{ij}))a_1 =  0$.\hfill $(Y_{1})$ \vspace{2mm}

\noindent
Here we have either $a_1=  0$ or
$$
\nu_{ij}^{(o,\pm)} =   -{n\over 2}\pm \sqrt{\({n\over 2}\)^2 + K
_{ij} -(n-1)}.
$$
Notice now that $\nu _{ij}^{e\pm}=  \nu^{o\pm }_{ij}+1$ and, more
then this, $K_{\nu^{o\pm }_{ij}+m+1}=  K_{\nu^{e\pm }_{ij}+m}$.
Comparing the power series solution $u_{ij}^{e\pm}(\ell )$ with
even indices $m$ and $u_{ij}^{o\pm}(\ell )$ with odd indices $m$
we see that $u_{ij}^{e\pm}(\ell )=  u_{ij}^{o\pm}(\ell )$.  Thus
it is enough to study the solution $u_{ij}^{e\pm}(\ell )$ only.
\vspace{2mm} \noindent We have for these solutions:
$$
\begin{array}{rcl}
\displaystyle \frac{n}{2}+\nu^{e\pm }_{ij}&= & \displaystyle 1\pm
\sqrt{\({n-2\over 2}\)^2 + K_{ij}}
\\
\\
&= &\displaystyle
1\pm
\frac{n-2}{2}\sqrt{1+\frac{\Lambda }{(n-1)(n-2)}
+\frac{8}{(n-2)^{2}}\left( \frac{\lambda ^{p}_{i}}{r_p^{2}}+
\frac{\lambda
^{q}_{j}}{r_q^{2}} \right) }
\\
\\
&= &\displaystyle
1\pm \frac{n-2}{2}\sqrt{\mu ^{2}+
\frac{8}{(n-2)^{2}}\left( \frac{\lambda ^{p}_{i}}{r_p^{2}}+
\frac{\lambda ^{q}_{j}}{r_q^{2}} \right) }.
\end{array}
$$
Here
$$
\mu ^{2}= 1+\frac{\Lambda
}{(n-1)(n-2)}= \frac{2}{(n-1)(n-2)}\[
\frac{p(p-1)}{r_p^{2}}+\frac{q(q-1)}{r_q^{2}} \]  ,
$$
see Appendix.
\vspace{2mm}

\noindent
From this we conclude that the solution $u^{e-}_{ij}$ with the leading
term $\ell ^{\nu _{ij}^{e-}}$ never belongs to $H^{1}_{2}(K)$ but
belongs to $L_{2}(K)$ if (and only if)
$$
\frac{n-2}{2}\sqrt{\mu ^{2}+
\frac{8}{(n-2)^{2}}\left( \frac{\lambda ^{p}_{i}}{r_p^{2}}+
\frac{\lambda ^{q}_{j}}{r_q^{2}} \right) }<1,
$$
which is equivalent to
$$
\mu ^{2}+
\frac{8}{(n-2)^{2}}\left( \frac{\lambda ^{p}_{i}}{r_p^{2}}+
\frac{\lambda ^{q}_{j}}{r_q^{2}} \right) <\frac{4}{(n-2)^{2}}.
$$
For $i= j= 0$ (radial solution) this condition, after substitution
of value for $\mu^2$, takes the form
$$
\left[ \frac{p(p-1)}{r_p^{2}}+\frac{q(q-1)}{r_q^{2}} \right]
<\frac{2(n-1)}{n-2}.
$$
These conditions on $n$, $p$, $r_p$, $r_q$ are met in the nonlinear
case as well; we call this situation the \emph{``minus-case''}.  Thus,
in the minus-case, the solution $u^{e-}_{00}$ belongs to
$L_{2}(K)$. On the other hand, the solution $u_{ij}^{e+}$ always
belongs to $H^{1}_{2}(K)$ (see Proposition 3.7 and Appendix). It
belongs to $H^{2}_{2}(K)$ (classical solution) if and only if
$$
\begin{array}{l}
\frac{n-2}{2}\sqrt{\mu ^{2}+
\frac{8}{(n-2)^{2}}\left( \frac{\lambda ^{p}_{i}}{r_p^{2}}+
\frac{\lambda ^{q}_{j}}{r_q^{2}} \right) }>1, \ \ \mbox{that is if}
\\
\\
\mu ^{2}+
\frac{8}{(n-2)^{2}}\left( \frac{\lambda ^{p}_{i}}{r_p^{2}}+
\frac{\lambda ^{q}_{j}}{r_q^{2}} \right) >\frac{4}{(n-2)^{2}}.
\end{array}
$$
For $i= j= 0$ (the radial solution) this condition, after substitution
of $\mu ^2$, takes the form
$$
\left[ \frac{p(p-1)}{r_p^{2}}+\frac{q(q-1)}{r_q^{2}} \right]
>\frac{2(n-1)}{n-2}.
$$
We will work with the same condition on $n,p,r_p,r_q$ in the nonlinear
case; we call this situation the \emph{``plus-case''}. Thus, in the
plus-case, the solution $u^{e+}_{00}$ belongs to $H^{2}_{2}(K)$.
\vspace{2mm}

\noindent
Now we return to the recursive equation $(Y_{m+2})$ for
the coefficients $a_{m}$ in the case $\nu _{ij}^{e+}$. We write
it in the form:
$$
K_{\nu^{e+}_{ij}+m}a_{m+2}=  -Q_{1}a_{m}, \ \ \ \mbox{which gives} \ \ 
\
a_{m+2}=  \frac{-Q_{1}a_{m}}{K_{\nu^{e+}_{ij}+m}}
$$
provided the denominator is not zero.
\vspace{2mm}

\noindent
Notice that, provided equation $(Y_{0})$ is satisfied, the expression
given above for $K_{\nu^{e+}_{ij}+m}$ can be rewritten as
$$
K_{\nu^{e+}_{ij}+m}= (m+2)(2\nu _{ij}+m+n).
$$
It follows from this formula that these coefficients are always
nonzero. We use the previous formula for $a_{m+2}$ recursively to
obtain
$$
a_{2m} =   {(-Q_1)^m\over \prod_{t=  1}^m
K_{\nu^{e+}_{ij}+2t}}a_0= {(-Q_1)^m\over \prod_{t=  1}^m
(2t+2)(2\nu_{ij}+2t+n)}a_0,
$$
and for the solution $u_{ij}^{e+}$, which we redenote to be
$u_{ij}$, we have
$$
u_{ij}= a_{ij}\ell ^{\nu
_{ij}}\sum_{m=  0}^{\infty} {(-Q_1)^m\over \prod_{t=  1}^m
(2t+2)(2\nu_{ij}+2t+n)}\ell ^{2m}
$$
with arbitrary constants $a_{ij}\in R$.
\vspace{2mm}

\noindent
We combine all calculations in the following theorem.
\begin{Theorem}\label{cone1-12a}
Let $M$ be a manifold with tame conical singularity as above, with
$\dim M \geq 5$, and $Q_1>0$.
\begin{enumerate}
\item[{\bf (1)}] There exists a solution $u_{ij}$ of the equation {\rm
(\ref{cone1-11})}, restricted on the cone $K$, is given by
\begin{equation}\label{cone1-13}
\begin{array}{l}
\displaystyle
u_{ij} =   \ell^{\nu_{ij}}\cdot \( \sum_{m=  0}^{\infty} {(-1)^m
Q_1^m \ell^{2m}\over \prod_{t= 1}^m (2t+2)(2\nu_{ij}+2t+n)}\)
\cdot a_0, \ \ \ \mbox{with}
\\
\\
\displaystyle
\nu_{ij}= \frac{n-2}{2}\[ \sqrt{
\mu^{2}+\frac{4}{(n-2)^{2}}\left(\frac{2\lambda
^{p}_{i}}{r_p^{2}+\frac{2\lambda ^{q}_{j}}{r_q^{2}}}\right)} -1
\] , \ \ \ \mbox{where}
\\
\\
\displaystyle
\mu ^{2}= \frac{2}{(n-1)(n-2)}\[
\frac{p(p-1)}{r_p^{2}}+\frac{q(q-1)}{r_q^{2}} \] .
\end{array}
\end{equation}
This solution belongs to the Sobolev space $H^{1}_{2}(M)$ and also
to the space $H^{2}_{2}(M)$ in the plus-case (see {\rm (\ref{plus})}).
\item[{\bf (2)}] The second linearly independent solution (denoted
above as $u_{ij}^{e,-}$) of the equation {\rm (\ref{cone1-11})} only 
belongs to $L_{2}(K)$ in the minus-case (see {\rm
(\ref{minus})}).
\item[{\bf (3)}] The general solution $u(\ell ,\theta, \psi )$ of
{\rm (\ref{cone1-11})} in $H^{1}_{2}(K)$ has the form
$$
u(\ell ,\theta, \psi )= \sum_{ij}a_{ij}\ell ^{\nu
^{e+}_{ij}}f_{ij}(\ell)\kappa ^{p}_{i}(\theta )\kappa^{q}_{j}(\psi ).
$$
Here the functions $f_{ij}(x)$ are defined by
$$
f_{ij}(\sqrt{Q_{1}}\ell )= \sum_{m=  0}^{\infty}{(-Q_1)^m\over 
\prod_{t=  1}^m
(2t+2)(2\nu_{ij}+2t+n)}\ell ^{2m},
$$
so that $u^{e+}_{ij}= a_{ij}\ell ^{\nu ^{e+}_{ij}}f_{ij}(\ell )$ and
the coefficients $a_{ij}$ ensure convergence of this series with
respect to $H^{1}_{2}(K)$-norm.
\end{enumerate}
\end{Theorem}
We notice that the exponent $\nu_{ij}$ of the solution $u_{ij}$ does not
depend on the eigenvalue $Q_1$.
\vspace{2mm}

\noindent
{\bf Remark.} We notice that the radial solution $u_{00}^{e+}(\ell )$ of 
the equation
(\ref{cone1-10}) has the following asymptotic
$$
u_{00}(\ell )\sim \ell ^{\nu _{00}},\ \nu _{00}= 
\frac{n-2}{2}\left( \sqrt{1+\frac{\Lambda}{(n-1)(n-2)}}-1\right) = \frac{n-2}{2}(\mu -1),
$$ and that for any solution $u(\ell ,\theta , \psi )$ of
(\ref{cone1-10}), there exists the radial solution $u_{0}(\ell )$ of
(\ref{cone1-10}) (radial part of $u$) such that as $\ell \rightarrow
0$,
$$
\vert u(\ell ,\theta , \psi )-u_{0}(\ell )\vert \leqq  C\ell
^{\sigma}, \ \ \ \mbox{for some $\sigma > 0$.}
$$
\section{Asymptotic of solutions: the nonlinear case}\label{s7}
In this section we study the nonlinear Yamabe equation near the
singular point:
\begin{equation}\label{cone1-30}
\frac{\p^2 u}{\p \ell^2}\! + \!{n-1\over \ell}{\p u \over \p\ell}\!+\! 
{2\over
{r_p^{2}\ell^2}}\Delta_{\theta}u\! +\!{2\over
{r_q^{2}\ell^2}}\Delta_{\psi}u\! -\! {n-2\over 4(n-1)} {\Lambda\over
\ell^2} u\! +\! Q_{\alpha}u^{\alpha}\!=\!0
\end{equation}
defined on the open set $(0,\epsilon)\times S^{p}\times S^{q}$.  Here
$3\leq p\leq n-3$, $\Lambda=\Lambda (p,q,r_p,r_q)$ is the curvature
factor, see Appendix for the details.
\vspace{2mm}

We restrict our attention to the \emph{radial solutions}, and return
to the general case at the end of this section.  Thus, we study
\emph{positive} solutions $u=u(\ell)$ of the equation
(\ref{cone1-30}).
\vspace{2mm}

\noindent
{\bf \ref{s7}.1. Reduction to a dynamical system.} We use the
cylindrical coordinates $t=-\ln \ell$, or $\ell = e^{-t}$, so that
$t\lra \infty$ as $\ell\lra 0$.  Then the equation (\ref{cone1-30})
becomes
\begin{equation}\label{cone1-39}
u_{tt} - (n-2)u_t - {n-2\over 4(n-1)} \Lambda u +
Q_{\alpha}e^{-2t}u^{\alpha}=0.
\end{equation}
We are looking for solutions of (\ref{cone1-39}) defined on the set
$$(-\ln(\epsilon),+\infty )\times S^{p}\times S^{q}.$$ We use the
substitution $u(t) = e^{\lambda t}w(t)$ in (\ref{cone1-39}), where
$\lambda ={2\over \alpha-1}$.  Then the time-dependence of the
coefficients in (\ref{cone1-39}) disappears, and one obtains the
following equation:
\begin{equation}\label{cone1-41}
\!\!\!\!\!\!\!\!\!
\begin{array}{l}
w^{\prime\prime}\! + \! \({4\over \alpha-1} \! - \! (n-2)\) 
w^{\prime}\!  +\! \({4\over
(\alpha-1)^2}\!  -\! {2(n-2)\over \alpha-1}\!  - \! 
{n-2\over 4(n-1)} \Lambda \) w\! 
+ \! Q_{\alpha}w^{\alpha}\! =\! 0.
\end{array}
\end{equation}
We exclude the case $\alpha=1$. Notice that the functions $u(\ell)$
and $w(t)$ are related as follows: $
u(\ell)=\ell^{-\frac{2}{\alpha-1}}w(-\ln\ell)$. We denote
\begin{equation}\label{cone1-41a}
\bar{b} = -{4\over \alpha-1} + (n-2), \ \ \ \bar{a}= -{4\over
(\alpha-1)^2} +{2(n-2)\over \alpha-1} + {n-2\over 4(n-1)} \Lambda .
\end{equation}
Let $x = w$, $y = w^{\prime}$. Then (\ref{cone1-41}) is equivalent to the
dynamical system
\begin{equation}\label{cone1-42}
\{\begin{array}{rcl}
x^{\prime} &=& y
\\
y^{\prime} &=& \bar{a} x + \bar{b} y - Q_{\alpha} x^{\alpha}.
\end{array}\right.
\end{equation}
{\bf \ref{s7}.2. The equilibrium points.} We find the equilibrium
points of (\ref{cone1-42}) by solving the system:
$$
\{\begin{array}{l}
y=0,
\\
\bar{a} x  - Q_{\alpha} x^{\alpha} = x(\bar{a} - Q_{\alpha} x^{\alpha-1})=0.
\end{array}\right.
$$
Since $y=0$, the second equation has the solution $x_1=0$ for all
values of the parameters and, in addition, the solution
\begin{equation}\label{cone1-41a-1}
x_2= \({\bar{a}\over Q_{\alpha}}\)^{{1\over \alpha-1}}>0
\end{equation}
provided ${\bar{a}\over Q_{\alpha}}>0$.  Thus the system
(\ref{cone1-42}) has one equilibrium point $w_1=(0,0)$ if
${\bar{a}\over Q_{\alpha}}\leq0$, and an additional one, $w_2=(x_2,0)$
with $x_{2}$ given by (\ref{cone1-41a-1}) if ${\bar{a}\over
Q_{\alpha}}>0$.
\vspace{2mm}

\noindent
{\bf Remark.} Notice that we consider only positive values for the root
$x_{2}$ since we are looking for the positive solutions of
(\ref{cone1-41}).  In terms of the dynamical system (\ref{cone1-42})
this means that a solution has to stay in the right half-plane for $t>T$
for some $T$.
\vspace{2mm}

\noindent
{\bf \ref{s7}.3. The equilibrium point $w_1$.} To study a behavior of
(\ref{cone1-42}) near $(0,0)$, we analyze its linear approximation:
$$
A_{(0,0)}= \(
\begin{array}{cc}
0 & 1 \\
\bar{a}& \bar{b}
\end{array}\)\ .
$$
The parameters ${\bar b}$ and ${\bar a}$ are defined in
(\ref{cone1-41a}). Thus, we have the characteristic equation
$$
\begin{array}{l}
\lambda ^2 - \bar{b}\lambda - \bar{a} = 0, \ \ \ \mbox{with the roots}
\ \ \ \
\lambda _{\pm} = {\bar{b}\over 2} \pm \sqrt{{\bar{b}^2\over 4}+ {\bar a}}.
\end{array}
$$
The corresponding eigenvectors $\vec{v}_{\pm}$ are given by
$$
\vec{v}_{\pm} = \(\begin{array}{c}
1 \\
\lambda _{\pm}
\end{array}\),  \ \ \mbox{so that} \ \ \ A_{(0,0)}\vec{v}_{\pm}=
\lambda _{\pm}\vec{v}_{\pm}.
$$
We show (see Appendix, Claim \ref{app1}) that $\bar{b}<0$ for all $1 <
\alpha \leq \alpha ^{*}$, and that ${\bar{b}^2\over 4}+ {\bar a} >0$
if $n\geq 3$. Therefore, the point $(0,0)$ is a {\sl saddle point} if
$\bar{a} >0$, and a {\sl stable focus} if $\bar{a}< 0$. \vspace{2mm}
\vspace{2mm}

\noindent
{\bf \ref{s7}.4. The equilibrium point $w_2$.}  Now we study the
second equilibrium point $w_2$:
$$
w_2 = (x_2,0), \ \ \ \mbox{with}\ \ \
x_2=\({\bar{a}\over
Q_{\alpha}}\)^{{1\over \alpha-1}},
$$
which exists provided $\bar{a}Q_{\alpha }>0$. Here we have the following
linear approximation of (\ref{cone1-42}) near the point $w_2$:
$$
A_{w_2}= \(
\begin{array}{cc}
0 & 1
\\
\bar{a}-\alpha Q_{\alpha}{\bar{a}\over
Q_{\alpha}}   & \bar{b}
\end{array}\) = \(
\begin{array}{cc}
0 & 1
\\
\bar{a}(1-\alpha) & \bar{b}
\end{array}\)
$$
with the characteristic equation $ r^2 - \bar{b}r + \bar{a}(\alpha-1)
=0 $.  We have the following eigenvalues $r_{\pm}$ and the
eigenvectors $\vec{z}_{\pm}$:
$$
\begin{array}{l}
r_{\pm} = {\bar{b}\over 2} \pm \sqrt{{\bar{b}^2\over 4}-
\bar{a}(\alpha-1)}, \ \ \ 
\vec{z}_{\pm} = \(\begin{array}{c}
1
\\
r_{\pm}
\end{array}\),  \ \ \mbox{so that} \ \ \ A_{(x_{2},0)}\vec{z}_{\pm}=
r_{\pm}\vec{z}_{\pm}.
\end{array}
$$
Thus, we have the following alternatives for the equilibrium point
$w_{2}$.
\begin{enumerate}
\item If ${\bar a}<0$, then $w_{2}$ is a saddle.
\item If ${\bar a}>0$, but $\frac{\bar{b}^2}{4}-
\bar{a}(\alpha-1)<0$, then the $w_{2}$ is a stable focus.
\item If ${\bar a}>0$, and ${\bar{b}^2\over 4}- \bar{a}(\alpha-1)>0$,
then $w_{2}$ is a stable node.
\end{enumerate}
It is shown in Appendix that all three cases are realized for different
values of $\alpha $.
\vspace{2mm}

\noindent
{\bf \ref{s7}.5. The phase pictures.} Now we study the critical points
of the system (\ref{cone1-42}) for different values of parameters
$\alpha$, $Q_{\alpha}$, and others. We determine asymptotic behavior
of the solutions of system (\ref{cone1-42}) by comparing them with the
corresponding solutions of the linearized system at the points
$w_{1},w_{2}$. For all $\alpha <\alpha ^{*}$ the points $w_{1},w_{2}$
are hyperbolic.  Thus, locally (near critical points) linear and
nonlinear phase pictures are trajectory-equivalent.  Moreover, those
solutions $w(t)$ of the nonlinear system which go to $w_{2}$ as $t
\rightarrow +\infty $ have the asymptotic behavior $w(t)\simeq w_{2}$,
and this determines to which Sobolev spaces they belong. For the
unbounded solutions ($w(t) \rightarrow \infty $), one has $x(t)
\rightarrow \infty $. Thus, the corresponding asymptotic behavior of
$u(\ell )$ is worse than that of an $\alpha$-basic function.  This
allows us to decide, in most cases, to which Sobolev space
$H^{k}_{2}(M)$ those solutions belong.
\vspace{2mm}

Finally, there are solutions $w(t)$ that tend to $w_{1}=(0,0)$ as $t
\rightarrow +\infty$. Locally (near $w_{1}$), the nonlinear system may
be thought of as a perturbation of the linear one. We use results by
Lettenmeyer, Hartman and Wintner (see \cite[Ch. 4, Theorems
5,9]{Coppel}, and \cite[Ch. X, Theorem 13.1, Corollary
16.3]{Hartman}). We check below that in our case the conditions of
those theorems are met.  These results guarantee that the principal
term of asymptotic behavior is the same for a solution tending to
the origin for the system (\ref{cone1-42}) and its linearization.
\vspace{2mm}

The results presented below depend on the inequalities below. We will
also describe these cases as the \emph{``plus-case''} and the
\emph{``minus-case''} respectively:
\begin{equation}\label{plus}
\begin{array}{l}
\mbox{{\bf The plus-case:}  \ \ \ \ \ \ \ $\displaystyle
\frac{p(p-1)}{r_p^{2}}+\frac{q(q-1)}{r_q^{2}}>\frac{2(n-1)}{n-2}$.}
\ \ \ \ \ \ \ \ \ \ \ \ \ \ 
\end{array}
\end{equation}
\begin{equation}\label{minus}
\!\!\!\!\!\!\!\!
\begin{array}{l}
\mbox{{\bf The minus-case:}  \ \ \ \ $\displaystyle
\frac{p(p-1)}{r_p^{2}}+\frac{q(q-1)}{r_q^{2}}<\frac{2(n-1)}{n-2} $.}
\ \ \ \ \ \ \ \ \ \ 
\end{array}
\end{equation}
We have the following alternative cases.
\vspace{2mm}

\noindent
{\bf Case 1:} ${\bar a}<0,\ Q_{\alpha }<0.$ In this case the phase
picture is given in Fig. \ref{s7}.1.  Here we have two families and
three separate solutions:
\vspace{2mm}

\noindent
{\bf (1)} The family $C_{\infty}$ consists of the solutions $w(t)$
going to $\infty $ asymptotically as $t \rightarrow +\infty $ nestling
at the unstable separatrix trajectory of the point $w_{2}$. We notice
that for given $\alpha $ the corresponding solution $u(\ell)$
approaches $+\infty$ faster, compared to an $\alpha$-basic function
$\ell^{-\frac{2}{\alpha-1}}$ as $\ell\to 0$.  
\vspace{2mm}

As it is proved in the Appendix, in the plus-case (\ref{plus}), $\alpha
_{0}<\frac{n+4}{n}$ and, since ${\bar a}<0$, $\alpha <\alpha
_{0}$. Thus, an $\alpha $-basic function does not belong to
$L_{2}(K)$. Therefore, the solutions from the family $C_{\infty
}$ do not belong to $L_{2}(K)$. 
\vspace{2mm}

In the minus case (\ref{minus}) if $\alpha <\frac{n+4}{n}$, the same
argument leads to the same conclusion, i.e. that the solutions from
the family $C_{\infty }$ do not belong to $L_{2}(K)$.  But if
$\frac{n+4}{n} <\alpha <\alpha _{0}$, then an $\alpha $-basic function
belongs to $L_{2}(K)$ and we do not know whether the solutions from
the family $C_{\infty }$ belong to $L_{2}(K)$.
\vspace*{32mm}

\begin{picture}(10,100)
\put(-51,-220){\PSbox{fig7-1.ps}{5cm}{100mm}}
\put(230,140){{\small $x^{\prime}\!\!=\!y$}}
\put(230,120){{\small $y^{\prime}\!\!=\!- 50 x\! -\! 15 y\! +\! 10 x^{1.2}$}}
\put(0,88){\vector(1,0){280}}
\put(290,88){{\small $x$}}
\put(33,10){\vector(0,1){170}}
\put(37,180){{\small $y$}}
\end{picture}
\vspace{2mm}

\centerline{{\small {\bf Fig. \ref{s7}.1.}
The phase picture for the case ${\bar a}<0$, $Q_{\alpha }<0$, here
$s=1/4$.}}
\vspace{2mm}

\noindent
{\bf (2)} The family $C_{0}$. Each solution of the family $w(t)\in
C_{0}$ goes to $w_{1}=(0,0)$ as $t \rightarrow +\infty $,
\emph{nestling at the direction of the eigenvector} $v_{+}$. Notice
that all but one (denoted by $w_{s}$, see below) of the admissible
solutions tending to $w_{1}$ belong to the family $C_{0}$. The
solutions $w(t)$ and the corresponding solutions $u(\ell )$ have the
following asymptotic behavior:
$$
\begin{array}{l}
w(t)\sim e^{\lambda _{+}t}= e^{\(\frac{\bar b}{2}+\sqrt{\frac{{\bar
b}^{2}}{4}+{\bar a}}\)t}, \ t\to \infty
\ \
u(\ell ) \sim \ell^{-\frac{2}{\alpha -1}-(\frac{\bar b}{2}+
\sqrt{\frac{{\bar b}^{2}}{4}}+\bar a)}, \ \ell\to 0.
\end{array}
$$
To analyze the behavior of solutions from the family $C_{0}$, we use
Lemma \ref{new-l1}.  We see that $u(\ell ) \sim \ell^{q}$ with $q=
-\frac{2}{\alpha -1}-(\frac{\bar b}{2}+\sqrt{\frac{{\bar
b}^{2}}{4}+{\bar a}})$.
\vspace{2mm}

Recall that $u(\ell )\in H^k_2(M)$ if and only if $k< \frac{n}{2}+q$.
Since ${\bar b}=(n-2)-\frac{4}{\alpha -1}$, we have $\frac{2}{\alpha
-1}+\frac{\bar b}{2}=\frac{n-2}{2}$; therefore, in this case
$$
\begin{array}{c}
\frac{n}{2}+q=\frac{n}{2}-\frac{n-2}{2}-\sqrt{\frac{{\bar
b}^{2}}{4}+{\bar a}}=1-\sqrt{\frac{{\bar b}^{2}}{4}+{\bar
a}}=1-\frac{n-2}{2}\mu , \ \ \mbox{where} \ \ 
\\
\\
\mu =\sqrt{1+\frac{\Lambda }{(n-1)(n-2)}} .
\end{array}
$$
As it is proved in the Appendix, such a solution never belongs to
$H^{1}_{2}(M)$.  It belongs to $L_{2}(M)$ only in the minus case
(\ref{minus}).
\vspace{2mm}

\noindent
{\bf (3)} Two incoming separatrix trajectories of the saddle point
$w_{2}$.  For both these trajectories, $w_{sep}(t) \rightarrow
w_{2}$. Thus, the corresponding solution $u_{sep}(\ell )$ is exactly
an $\alpha $-basic function near $x^*$. If in this case $\alpha
<\alpha _{0}<\frac{n+4}{n}$ (in the plus-case (\ref{plus}), see the
Appendix), the solutions $u_{sep}(\ell )$ do not belong to $L_{2}(M)$
(recall that here $\alpha <\alpha _{0}$ since ${\bar a}<0$).
\vspace{2mm}

\noindent
{\bf (4)} The solution $w_{s}(t)\rightarrow w_{1}=(0,0)$,
corresponding to the eigenvector $v_{-}.$ This solution, obtained by a
$C^{1}$-diffeomorphic twist of the corresponding solution of the
linearization of (\ref{cone1-41a-1}) at $w_{1}$, has the asymptotic
behavior $w_{s}(t) \sim e^{\lambda _{-}t}$. The corresponding solution
$u_{s}(\ell )\sim l^{q},\ q=-\frac{2}{\alpha -1}-\lambda _{-}$, where
$$
\begin{array}{c}
q=-\frac{2}{\alpha -1}-\frac{\bar b}{2}+\sqrt{\frac{{\bar b}^{2}}{4}+{\bar a}}
=-\frac{n-2}{2}+\sqrt{\frac{{\bar b}^{2}}{4}+{\bar a}},
\end{array}
$$
see above. Therefore, 
$$
\begin{array}{c}
\frac{n}{2}+q=1+\sqrt{\frac{{\bar
b}^{2}}{4}+{\bar a}}=1+\frac{n-2}{2}\mu .
\end{array}
$$ 
As it is shown in
the Appendix, the solution $u_{s}(t)$ always belongs to $H^{1}_{2}(M)$.
Furthermore, $u_{s}(\ell )\in H^{2}_{2}(M)$ if and only if the plus
case (\ref{plus}) condition is met.
\vspace{2mm}

\noindent
{\bf Case 2:} ${\bar a}<0,\ Q_{\alpha }>0.$ In this case the phase
picture is given in Fig. \ref{s7}.2.
Here we have one family $C_{0}$ of admissible solutions and a special
solution $u_{s}(\ell )$.  Solutions of the family $C_{0}$ tend to
$w_{1}=(0,0)$ as $t \to +\infty$. Such solutions $w(t)$ asymptotically
nestle in the direction of the eigenvector $v_{+}$. 
\vspace{2mm}

The asymptotic
behavior of $w(t)$ and of the corresponding solutions $u(\ell)$ is the
following:
$$
\begin{array}{l}
w(t)\sim e^{\lambda _{+}t}= e^{\(\frac{\bar b}{2}+\sqrt{\frac{{\bar
b}^{2}}{4}+{\bar a}}\)t}, \ t\to \infty,
\ \
u(\ell ) \sim \ell^{-\frac{2}{\alpha -1}-(\frac{\bar b}{2}+
\sqrt{\frac{{\bar b}^{2}}{4}}+\bar a)}, \ \ell \to 0.
\end{array}
$$
The same argument as in Case 1 proves that solutions from the
family $C_{0}$ belong to $L_{2}(M)$ only in the minus-case
(\ref{minus}) and that none of these solutions belong to $H^{1}_{2}(M)$.
Similar to Case 1, the solution $u_{s}(\ell )$ always belongs to
$H^{1}_{2}(M)$ and belongs to $H^{2}_{2}(M)$ if and only if the
plus-case condition (\ref{plus}) is met.
\vspace*{31mm}

\begin{picture}(10,100)
\put(-51,-220){\PSbox{fig7-2.ps}{5cm}{100mm}}
\put(230,140){{\small $x^{\prime}\!\!=\!\!y$}}
\put(230,120){{\small $y^{\prime}\!\!=\!\!- 50 x\! - \!15 y \!- \!10 x^{1.2}$}}
\put(0,88){\vector(1,0){280}}
\put(290,88){{\small $x$}}
\put(33,10){\vector(0,1){170}}
\put(37,180){{\small $y$}}
\end{picture}
\vspace{2mm}

\centerline{{\small {\bf Fig. \ref{s7}.2.}  The phase picture for the
case ${\bar a}<0,\ Q_{\alpha }>0$, here $s=1/4$.}}
\vspace{2mm}

\noindent
{\bf Case 3: ${\bar a}>0,\ Q_{\alpha }>0,\ \frac{{\bar
b}^{2}}{4}-{\bar a}(\alpha -1)<0.$} In this case the phase
picture is given in Fig. \ref{s7}.3.
\vspace*{31mm}

\begin{picture}(10,95)
\put(-51,-220){\PSbox{fig7-3.ps}{5cm}{100mm}}
\put(232,140){{\small $x^{\prime}\!\!=\!\!y$}}
\put(232,120){{\small $y^{\prime}\!\!=\!\!5.5 x\! - \!(5/3) y\! -
\! 5 x^{1.6}$}}
\put(0,88){\vector(1,0){280}}
\put(290,88){{\small $x$}}
\put(33,10){\vector(0,1){170}}
\put(37,180){{\small $y$}}
\end{picture}
\vspace{2mm}

\centerline{{\small {\bf Fig. \ref{s7}.3.}  The phase picture for the
case ${\bar a}>0,\ Q_{\alpha }>0$,}}
\centerline{{\small $\frac{{\bar b}^{2}}{4}-{\bar
a}(\alpha -1)<0$, here $s=1/4$.}}
\vspace{2mm}

\noindent
\noindent
{\bf Case 3$^{\prime}$: ${\bar a}>0,\ Q_{\alpha }>0,\ \frac{{\bar
b}^{2}}{4}-{\bar a}(\alpha -1)>0.$} We have the following phase given at
Fig. \ref{s7}.3$^{\prime}$.
\vspace*{32mm}

\begin{picture}(10,100)
\put(-45,-220){\PSbox{fig7-3a.ps}{5cm}{100mm}}
\put(237,140){{\small $x^{\prime}\!\!=\!\!y$}}
\put(237,120){{\small $y^{\prime}\!\!=\!\!5.5 x\!\! - \!\!2.15 y\!\! -\!\! 5 x^{1.56}$}}
\put(0,88){\vector(1,0){280}}
\put(290,88){{\small $x$}}
\put(33,10){\vector(0,1){170}}
\put(37,180){{\small $y$}}
\end{picture}
\vspace{2mm}

\centerline{{\small {\bf Fig. \ref{s7}.3$^{\prime}$.} The phase
picture for the case ${\bar a}>0,\ Q_{\alpha }>0,\ \frac{{\bar
b}^{2}}{4}-{\bar a}(\alpha -1)>0.$}}
\vspace{2mm}

\noindent
We analyze the cases 3, 3$^{\prime}$ together since they have very similar
classes of admissible solutions.
\vspace{2mm}

Here we have the class $C_{F}$ of Fowler solutions (see
\cite{Bateman}): Here $w(t)\to w_{2}$ as $t \to +\infty $, and
$$
u(\ell )\sim \ell ^{-\frac{2}{\alpha -1}}, \ \ \ \ell\to 0.
$$
Also we have a separatrix solution $w_{s}(t)$ that approaches
$w_{1}=(0,0)$ as $t\to +\infty $.  This solution asymptotically
nestles at the eigendirection $v_{-}$. It has the following
asymptotic:
$$
\begin{array}{l}
w_{s}(t)\!\sim \!e^{\lambda _{-}t}\!= \!e^{\(\frac{\bar b}{2}-\sqrt{\frac{\bar
{b}^{2}}{4}+{\bar a}}\)t}\!\!, \ t\to+\infty,
\ \
u_{s}(\ell )\! \sim \!\ell^{-\frac{2}{\alpha -1}-(\frac{\bar
b}{2}-\sqrt{\frac{{\bar b}^{2}}{4}+{\bar a}})}, \ \ell\to 0.
\end{array}
$$

\noindent
{\bf Remark.} In fact, the solution $u_{s}(\ell )$ gives a minimum for
the corresponding Yamabe functional. This asymptotic behavior refines the
general result on the Yamabe minimizer for cylindrical manifolds for
the particular ``slice'' $S^p\times S^q$, see \cite{AB}.
\vspace{2mm}

\noindent
The Fowler solutions in both cases 3 and 3$^{\prime}$ do not belong to
$H^{1}_{2}(M)$, by Proposition \ref{new-p1}. Moreover, they even do not
belong to $L_{2}(M)$ for $\alpha _{0} < \alpha \leq
\frac{n+4}{n}$ (the minus-case (\ref{minus})); however, they belong to
$L_{2}(M)$ for $\frac{n+4}{n}<\alpha $.
\vspace{2mm}

The separatrix solution $u_{s}(\ell )$ has the asymptotical behavior
$\ell ^{q}$ as $\ell \rightarrow +\infty $, with $ q=-\frac{2}{\alpha
-1}-\({\frac{\bar b}{2}-\sqrt{\frac{{\bar b}^{2}}{4}+{\bar a}}}\).
$
Using $\frac{2}{\alpha
-1}+\frac{\bar b}{2}=\frac{n-2}{2}$, we obtain
$$
\begin{array}{l}
q=
-\frac{n-2}{2}+\sqrt{\frac{{\bar b}^{2}}{4}+{\bar
a}}=-\frac{n-2}{2}+\frac{n-2}{2}\mu.
\end{array}
$$
Thus, $\frac {n}{2}+q=1+\frac{n-2}{2}\mu$.  Similarly to the Case 1,
solution $u_{s}(\ell )$ always belongs to $H^{1}_{2}(M)$ but belongs
to $H^{2}_{2}(M)$ if and only if the plus-case condition (\ref{plus})
is met.
\vspace*{33mm}

\begin{picture}(10,100)
\put(-51,-220){\PSbox{fig7-4.ps}{5cm}{100mm}}
\put(231,140){{\small $x^{\prime}\!\!=\!\!y$}}
\put(231,120){{\small $y^{\prime}\!\!=\!\!5.5 x \!\!- \!\!\(5/3\) 
y\!\! +\!\! 5 x^{1.6}$}}
\put(0,88){\vector(1,0){280}}
\put(290,88){{\small $x$}}
\put(33,10){\vector(0,1){170}}
\put(37,180){{\small $y$}}
\end{picture}
\vspace{2mm}

\centerline{{\small {\bf Fig. \ref{s7}.4.} The phase picture for the
case ${\bar a}>0,\ Q_{\alpha }<0$, here $s=3/4$.}}
\vspace{2mm}

\noindent
{\bf Case 4: ${\bar a}>0,\ Q_{\alpha }<0.$} Here the phase picture is
given in Fig. \ref{s7}.4.  In this case we have the family $C_{\infty
}$ of admissible solutions $w(t)$, such that $w(t)\to\infty$ as
$t\to\infty $.  The corresponding solution $u(\ell )$ goes to $\infty$
faster then the \emph{$\alpha$-basic function} $\ell^{-\frac{2}{\alpha
-1}}$ for a given $\alpha$.  Thus, these solutions do not belong to
$H_{2}^{1}(M).$
\vspace{2mm}

Here, $\alpha >\alpha _{0}$, and in the plus case (\ref{plus}), or in
the minus case (\ref{minus}) when $\alpha >\frac{n+4}{n}$, an $\alpha
$-basic function belongs to $L_{2}(K)$. But, we do not know if
solutions from $C_{\infty}$ belong to this space. On the other hand,
in the minus-case (\ref{minus}), if $\alpha <\frac{n+4}{n}$, such
solutions do not belong to $L_{2}(K)$.

\vspace{2mm}

\noindent
In addition, we have the separatrix solution $w_{s}(t)$ that tends to
$w_{1}=(0,0)$ as $t \to+\infty$. This solution asymptotically nestles
in the direction of the eigenvector $v_{-}$. The asymptotic of
$w_{s}(t)$ and $u_{s}(\ell )$ is the following:
$$
\begin{array}{l}
w_{s}(t)\!\sim\! e^{\lambda _{-}t}= e^{\(\frac{\bar b}{2}-\sqrt{\frac{\bar
{b}^{2}}{4}+{\bar a}}\)t}, \ t\to\infty,
\ \
u_{s}(\ell ) \!\sim\! \ell ^{-\frac{2}{\alpha -1}-(\frac{\bar
b}{2}-\sqrt{\frac{{\bar b}^{2}}{4}+{\bar a}})}, \ \ell \to 0.
\end{array}
$$
Thus, the separatrix solution $u_{s}(\ell )$ has asymptotic $\ell ^{q}$ with
$$
q=-\frac{2}{\alpha -1}-({\frac{\bar b}{2}-\sqrt{\frac{{\bar
b}^{2}}{4}+{\bar a}}})=-\frac{n-2}{2}+\sqrt{\frac{{\bar b}^{2}}{4}+{\bar a}},
$$
and $\frac{n}{2}+q=1+\frac{n-2}{2}\mu.$ Similarly to the Case 1,
solution $u_{s}(\ell )$ always belongs to $H^{1}_{2}(M)$ but belongs
to $H^{2}_{2}(M)$ if and only if the plus-case condition (\ref{plus})
is met.
\vspace{2mm}

Now we consider the special cases that were previously left outside of
the scope.
\vspace{2mm}

\noindent
{\bf Case 5.  $Q_{\alpha}=0$.}  In this case, we have a linear system
with the matrix
$$
A=\left( \begin{array}{cc} 0 & 1 \\ {\bar a} &
{\bar b}
\end{array} \right) .
$$
Thus, the solutions can be written out explicitly.  Depending on the
sign of ${\bar a}$, we have either a stable node (if ${\bar a} <0$,
Case 5$-$, Fig. \ref{s7}.5$-$) or a saddle point (if ${\bar a} >0$,
Case 5$+$, Fig. \ref{s7}.5$+$).  Correspondingly, we have linear
versions of cases 2 and 4. Conclusions that were made in these cases
regarding asymptotics of admissible solutions are true in this case
($Q_{\alpha }=0$) as well.

In the Case $5+$, where ${\bar a}>0$, the point $w_{1}$ is a saddle.
A generic admissible solution from the family $C_{\infty }$ is growing
faster than an $\alpha$-basic function.  More specifically, we have
$w(t) \sim e^{\lambda _{+}t}$ as $t \to +\infty $. Thus, the
corresponding solution $u(\ell )$ has the asymptotic:
$$
u(\ell )\sim \ell ^{-\frac{2}{\alpha -1}-\lambda _{+}}, \ \ \ \ell\to 0.
$$
We use Lemma \ref{new-l1}: we have $q= -\frac{2}{\alpha
-1}-(\frac{\bar b}{2}+ \sqrt{\frac{{\bar b}^{2}}{4}+{\bar a}})$.  The
same argument as before gives that
$
\frac{n}{2}+q=1-\sqrt{\frac{\bar{b}^{2}}{4}+{\bar
a}}=1-\frac{2}{n-2}\mu .
$ 
Thus, this solution never belongs to
$H^{1}_{2}(K)$, and it does belong to $L_{2}(K)$ if and only if the
minus-case condition (\ref{minus}) is met.
\vspace*{31mm}

\begin{picture}(10,100)
\put(-50,-220){\PSbox{fig7-5a.ps}{5cm}{100mm}}
\put(240,140){{\small $x^{\prime}=y$}}
\put(240,120){{\small $y^{\prime}=5.5 x - (5/3) y$}}
\put(0,88){\vector(1,0){280}}
\put(290,88){{\small $x$}}
\put(33,10){\vector(0,1){170}}
\put(37,180){{\small $y$}}
\end{picture}
\vspace{2mm}

\centerline{{\small {\bf Fig. \ref{s7}.5$+$.} The phase picture for
the case ${\bar a}>0,\ Q_{\alpha }=0$, here $s=3/4$.}}

\vspace*{33mm}

\begin{picture}(10,100)
\put(-50,-220){\PSbox{fig7-5b.ps}{5cm}{100mm}}
\put(240,140){{\small $x^{\prime}=y$}}
\put(240,120){{\small $y^{\prime}=-50 x - 15 y$}}
\put(0,88){\vector(1,0){280}}
\put(290,88){{\small $x$}}
\put(33,10){\vector(0,1){170}}
\put(37,180){{\small $y$}}
\end{picture}
\vspace{2mm}

\centerline{{\small {\bf Fig. \ref{s7}.5$-$.} The phase picture for
the case ${\bar a}<0,\ Q_{\alpha }=0$, here $s=1/4$.}}
\vspace{2mm}

\noindent
\noindent
We also have the stable separatrix solution $u_{s}(\ell )$, which has
the asymptotic $\ell ^{q}$ with
$$
q=-\frac{2}{\alpha -1}-\left(\frac{\bar b}{2}-\sqrt{\frac{{\bar
b}^{2}}{4}+{\bar a}}\right).
$$
Similarly to the above cases, we conclude that $u_{s}(\ell )\in
H^{1}_{2}(M)$, and it also belongs to $H^{2}_{2}(M)$ if and only if
the plus-case condition (\ref{plus}) is met.
\vspace{2mm}

In the Case $5-$ (where ${\bar a}<0$) we have the linear version of
the Case 2.  Here, the generic (slow) admissible solutions from the
family $C_{0}$, which are nestling in the eigendirection $v_{+}$, never
belong to $H^{1}_{2}(M)$ and belong to $L_{2}(M)$ if and only if the
minus-case condition (\ref{minus}) is met.
\vspace{2mm}

The fast-decreasing solution $u_{s}(t)=\ell ^{{-\frac{2}{\alpha
-1}}-\lambda _{-}},$ corresponding to the smaller eigenvalue, has the
asymptotic $\sim \ell ^{1+\sqrt{\frac{{\bar b}^{2}}{4}+{\bar a}}}=\ell
^{1+\frac{2}{n-2}\mu }$.  Thus, according to the discussion in the
Appendix, it always belongs to $H^{1}_{2}(M)$ and belongs to
$H^{2}_{2}(M)$ if and only if the plus-case condition (\ref{plus}) is
met.
\vspace{2mm}

\noindent
{\bf Case 6: ${\bar a}=0.$} In this case, (when $\alpha =\alpha _{0}$,
see Appendix) $w_{1}$ is the only singular point of the system
(\ref{cone1-42}). This point is degenerate: $ \lambda _{-}={\bar
b}_{0}=-(n-2)\mu$, $\lambda _{+}=0$. Depending on the sign of
$Q_{\alpha }$ we have the cases $6+$ or $6-$.
\vspace{2mm}

\noindent
{\bf Case 6$+$: ($Q_{\alpha }>0$, a weak stable node).} In this case
the phase picture is the following:
\vspace*{33mm}

\begin{picture}(10,100)
\put(-45,-220){\PSbox{fig7-6a.ps}{5cm}{100mm}}
\put(250,140){{\small $x^{\prime}=y$}}
\put(250,120){{\small $y^{\prime}=-5y -5x^{1.4}$}}
\put(0,88){\vector(1,0){280}}
\put(290,88){{\small $x$}}
\put(33,10){\vector(0,1){170}}
\put(37,180){{\small $y$}}
\end{picture}
\vspace{2mm}

\centerline{{\small {\bf Fig. \ref{s7}.6+.} The phase picture for
$\bar{a}=0$, $Q_{\alpha}>0$, here $s=1/2$.}}
\vspace{2mm}

\noindent
This case is similar to Case 2. The difference is in the fact that
$\lambda _{+}=0.$ As a result, the solutions of the family $C_{0}$
tend to $w_{1}=(0,0)$ slowly and, similar to the Case 2, do not belong
to $L_{2}(M)$.
\vspace{2mm}

The fast solution $w_{s}(t)$ corresponds to the second eigenvalue
$\lambda _{-}={\bar b}$, where the corresponding eigenfunction
$u_{s}(\ell )$ behaves asymptotically as $u_{s}(\ell )\sim l^{q}$ ,
with $q=-\frac{2}{\alpha_{0}-1}-{\bar b}$. Thus, we have: $
\frac{n}{2}+q = \frac{n}{2}-\frac{2}{\alpha _{0}-1}-{\bar b}
=1-\frac{\bar b}{2}=2-\frac{n}{2}+\frac{2}{\alpha
_{0}-1}=1+\frac{n-2}{2}\mu.  $ Therefore, $u_{s}(\ell )\in
H^{1}_{2}(M)$ and, in addition, it belongs to $H^{2}_{2}(M)$ if the
plus-case condition (\ref{plus}) is met.
\vspace{2mm}

\noindent
{\bf Case 6-: ($Q_{\alpha }<0$, weak saddle).} In this case the phase
picture is given at Fig. \ref{s7}.6-.  The Case 6$-$ is similar to the
Case 4. We have the family of solutions $C_{\infty}$ and a separatrix
solution $w_s(t)$.
\vspace{2mm}

The family $C_{\infty }$ consists of solutions $u(\ell )$ which grow
faster than the basic solutions $\ell^{-\frac{2}{\alpha _{0}
-1}}$. Thus, these solutions do not belong to $L_{2}(M)$ if the value
$\alpha _{0}$ (for which ${\bar a}=0$) is smaller then the critical
value $\frac{n+4}{n}$ from Proposition \ref{new-p1}, that is, in the
minus-case (\ref{minus}). In the opposite case, we do not know if
solutions from the family $C_{\infty}$ belong to $L_{2}(K)$ or not.
\vspace{2mm}

On the other hand, separatrix solution $u_{s}(\ell )$ has the
asymptotic $\ell ^{q}$ with $q=-\frac{2}{\alpha -1}-{\bar b}.$
Similarly to the above consideration we conclude that separatrix
solution $u_{s}(\ell )\in H^{1}_{2}(M)$ and it belongs to
$H^{2}_{2}(K)$ if the plus-case condition (\ref{plus}) is met.
\vspace{2mm}
\vspace*{33mm}

\begin{picture}(10,100)
\put(-45,-220){\PSbox{fig7-6b.ps}{5cm}{100mm}}
\put(250,140){{\small $x^{\prime}=y$}}
\put(250,120){{\small $y^{\prime}= - 5 y + 5 x^{1.4}$}}
\put(0,88){\vector(1,0){280}}
\put(290,88){{\small $x$}}
\put(33,10){\vector(0,1){170}}
\put(37,180){{\small $y$}}
\end{picture}
\vspace{2mm}

\centerline{{\small {\bf Fig. \ref{s7}.6-.} The phase picture for
$\bar{a}=0$, $Q_{\alpha}<0$, here $s=1/2$.}}
\vspace{2mm}

\noindent
{\bf Case 7. The critical: $\alpha =\alpha ^{*}=\frac{n+2}{n-2}.$} In
this case $\alpha -1=\frac{4}{n-2},\ s=1.$ We have here ${\bar b}=0$,
and ${\bar a}(1)=\bar{a}^*=
\frac{(n-2)^{2}}{4}(1+\frac{\Lambda}{(n-2)(n-1)})
=\frac{(n-2)^{2}}{4}\mu ^{2}$. Therefore (see Appendix)
$\bar{a}^*>0$. Our system takes the form
\begin{equation}\label{cone1-4?}
\{\begin{array}{rcl}
x^{\prime} &=& y
\\
y^{\prime} &=& {\bar{a}}^{*} x - Q_{{\alpha}^{*}} x^{{\alpha}^{*}}.
\end{array}\right.
\end{equation}
We have the corresponding equation of the second order
$$
w''-\bar{a}^*w+Q_{\alpha ^*}w^{\alpha ^*}=0.
$$
This equation has a first integral (a Hamiltonian function of system
(\ref{cone1-4?})):
$$
I(x,y)=\frac{y^2}{2}-\frac{\bar{a}^*x^{2}}{2}+\frac{Q_{{\alpha}^{*}}}{{\alpha}^{*}+1}
x^{{\alpha}^{*}+1}
$$
For the singular point $w_{1}=(0,0)$ we have $\lambda _{\pm }=\pm
\sqrt{{\bar a}^{*}}$, thus $w_{1}$ is saddle point.
\vspace{2mm}

\noindent
{\bf Case 7$-$: $Q_{{\alpha}^{*}}<0$.} In this case $w_{1}$ is the
only singular point.  The phase picture is given at Fig. \ref{s7}.7-.

\vspace*{33mm}

\begin{picture}(10,100)
\put(-40,-220){\PSbox{fig7-7a.ps}{5cm}{100mm}}
\put(250,140){{\small $x^{\prime}=y$}}
\put(250,120){{\small $y^{\prime}=x+ x^{1.8}$}}
\put(0,88){\vector(1,0){280}}
\put(290,88){{\small $x$}}
\put(33,10){\vector(0,1){170}}
\put(37,180){{\small $y$}}
\end{picture}
\vspace{2mm}

\centerline{{\small {\bf Fig. \ref{s7}.7-.} The phase picture for the
case $\alpha =\alpha ^{*},\ Q_{\alpha ^*}<0.$}}
\vspace{2mm}

\noindent
Here we have the family $C_{\infty }$ of admissible solutions and the
incoming separatrix $w_{s}(t)$. All solutions of the family
$C_{\infty}$ go to infinity as $t \to +\infty $. For such solutions
$w(t) \to +\infty $ as $t \to +\infty $, and $u(\ell )$ goes to $+
\infty $ faster than the basic function $\ell^{-\frac{2}{\alpha
^{*}-1}}$ as $\ell \to 0$. Thus, none of the solutions $u(\ell )\in
C_{\infty }$ belong to $H_{2}^{1}(M)$. Furthermore, since $\alpha
^{*} > \frac{n+4}{4}$, an $\alpha ^{*}$-basic function does belong to
$L_{2}(K)$, and in order to determine whether $u(\ell )$ belongs to
$L_2(M)$, we have to study its asymptotic behavior in more detail.

The separatrix solution $w_{s}(t)$ that tends to
$w_{1}=(0,0)$ as $t \to +\infty .$ This solution
asymptotically nestles in the direction of the eigenvector $v_{-}$. We have the
asymptotic:
$$
\begin{array}{l}
w_{s}(t)\sim e^{\lambda _{-}t}= e^{-\sqrt{\bar{a}^*}t}, \ \ \ t\to\infty,
\ \
\ \
u_{s}(\ell ) \sim \ell ^{-\frac{n-2}{2}+\sqrt{\bar{a}^*}}\ \ \ \ell\to0.
\end{array}
$$
For this solution $\frac{n}{2}+q=1+\frac{n-2}{2}\mu$. Thus, this
solution always belongs to $H^{1}_{2}(M)$ and does belong to $H_2^2(M)$
if the plus-case condition (\ref{plus}) is fulfilled.

\vspace{2mm}

\noindent
{\bf Case 7$+$: $Q_{{\alpha}^{*}}>0$.}
If $Q_{\alpha ^{*}}>0$, we have the second singular point
$w_{2}=x_{2}=((\frac{\bar{a}^*}{Q_{\alpha ^{*}}})^{\frac{n-2}{4}} ,0).$ The
eigenvalues of the linearization at this point are
$$
r_{\pm }=\pm \sqrt{-\frac{4\bar
a^{*}}{n-2}}.
$$
Thus, $w_{2}$ is a center, see Fig. \ref{s7}.7+ for the phase picture.
\vspace*{33mm}

\begin{picture}(10,100)
\put(-40,-220){\PSbox{fig7-7b.ps}{5cm}{100mm}}
\put(250,140){{\small $x^{\prime}=y$}}
\put(250,120){{\small $y^{\prime}=- 5 y - 5 x^{1.6}$}}
\put(0,88){\vector(1,0){280}}
\put(290,88){{\small $x$}}
\put(33,10){\vector(0,1){170}}
\put(37,180){{\small $y$}}
\end{picture}
\vspace{2mm}

\centerline{{\small {\bf Fig. \ref{s7}.7+.} The phase picture for the
case $\alpha =\alpha ^{*},\ Q_{\alpha ^*}>0.$}}

\vspace{2mm}

\noindent
Thus, inside the homoclinic loop we have a family of periodic
solutions. These solutions are known as Fowler or Delaunay solutions
(see \cite{KMPS}). Notice that the value of the first integral $I$ on
the separatrix loop is zero, and, at $w_{2}$,
$$
I(w_{2})=-\frac{Q_{\alpha
^{*}}}{n}\left(\frac{\bar{a}^*}{Q_{\alpha ^{*}}}\right)
^{\frac{n}{2}}<0.
$$
Therefore, a Fowler solution is determined by the value of the
integral $I$ in the interval $(I(w_{2}),0)$, or by the minimal value
$w_{\min}$ of $w(t)=(x(t),y(t))$ along its trajectory.
The integral $I(x,y)$ takes values in $[I(w_2),0]$. Such a solution is
also determined by the minimum value $x_{\min}$ of $x(t)$,
$x_{\min}\in (0,(\frac{{\bar a}^{*}}{Q_{\alpha
^{*}}})^{\frac{n-2}{4}})$.  For such a trajectory $x_{\min}\leq x(t)
\leq x_{\max} $ for all $t$.  As a result, for the corresponding
function $u(\ell) $ we have
$$
x_{\min}\ell ^{-\frac{n-2}{2}} \leq u(\ell ) \leq x_{\max}\ell
^{-\frac{n-2}{2}} .
$$
Therefore, the Fowler solutions have the asymptotic
$\ell^{-\frac{n-2}{2}}$ as $\ell \rightarrow 0$ and (see Proposition
\ref{new-p1}) do belong to $L_{2}(M)$, but not to $H^{1}_{2}(M)$.

In addition, we have the incoming separatrix solution $w_{s}(t)$ that
tends to $w_{1}=(0,0)$ as $t \rightarrow +\infty$. This solution
asymptotically nestles in the direction of the eigenvector $v_{-}$. We
have the asymptotic
$$
\begin{array}{l}
w_{s}(t)\sim e^{\lambda _{-}t}= e^{-\sqrt{\bar{a}^*}t}, \ \ \ t\to\infty
\ \
\ \
u_{s}(\ell ) \sim \ell ^{-\frac{n-2}{2}+\sqrt{\bar{a}^*}}, \ \ \ \ell\to 0.
\end{array}
$$
The separatrix solution $u_{s}(\ell )\sim \ell ^{q}$ with $
q=-\frac{n-2}{2}+\sqrt{{\bar a}^*}=\frac{n-2}{2}(\mu -1).  $ Thus,
$\frac{n}{2}+q=1+\frac{n-2}{2}\mu $, and this solution always
belongs to $H^{1}_{2}(M)$ and belongs to $H^{2}_{2}(M)$ if and only if
the plus-case condition (\ref{plus}) is fulfilled.
\vspace{2mm}

We denote
\begin{equation}\label{ppp1}
\displaystyle \sigma
=\frac{n-2}{2}(\mu -1),\ \ \ \ \mu ^{2}=
\frac{2}{(n-1)(n-2)}\left[\frac{p(p-1)}{r_p^{2}}+\frac{q(q-1)}{r_q^{2}}\right].
\end{equation}
{\bf \ref{s7}.6. Perturbation near the point $w_1$.}  Now we consider
the system (\ref{cone1-42}) as a perturbation of its linearization
near $w_{1}$.  We will apply \cite[Ch. 4, Theorems 5 and 9]{Coppel} or
\cite[Ch. X, Corollary 16.3]{Hartman} to those solutions, which tend to
zero as $t\to \infty$. We claim that two following conditions for the
nonlinear term $Q_{\alpha}x^{\alpha }$ (with $\alpha >1$) are
satisfied.
\vspace{2mm}

\noindent
{\bf (1)} Indeed, the perturbation $f(t,x)$ must satisfy
$$
\vert f(t,x) \leq L \vert x\vert^{1+\rho }
$$
for $t\geq t_{0},\vert x\vert \leq \delta $, where $L,\delta ,\rho $
are positive constants. In our case, it is enough to take
$\rho =\alpha -1$.
\vspace{2mm}

\noindent
{\bf (2)} The Lipschitz condition for $f(t,x)$ is satisfied since
$x^{\alpha }$ is differentiable for $x\geq 0$ and its derivative
$\alpha x^{\alpha -1}$ can be made arbitrarily small for $\vert x \vert
<\delta $ with small enough $\delta $.
\vspace{2mm}

Thus \cite[Ch. 4, Theorems 5 and 9]{Coppel} imply that there is
one-to-one correspondence between those solutions of the system
(\ref{cone1-42}) that tend to zero as $t\rightarrow +\infty $ and of
the corresponding solution of the linearization of
(\ref{cone1-42}). This shows that the norms $\vert w(t)\vert$ of
corresponding solutions have the same asymptotical behavior as
$t\rightarrow +\infty $. Using more refined results of (Hartman,
Theorem 13.1 and Corollary 16.2) we may conclude the same about the
asymptotic behavior of solutions $w(t)$ themselves, rather then their
norms).
\vspace{2mm}

We have proved the following results.
\begin{Theorem}\label{nonl}
Let $M$ be a manifold with tame conical singularities as above, $n\geq
5$, $K\subset M$ be its conical part, and $1\leq \alpha \leq
\alpha^*$.
\begin{enumerate}
\item[{\bf (1)}] Then the equation {\rm (\ref{cone1-39})} has unique
radial solution $u_s(\ell)$ which belongs to the space
$H_{2}^{1}(K)$. Moreover, $ u_s(\ell) \sim \ell^{-\sigma}$, where
$\sigma =\frac{n-2}{2}(\mu -1)$ (see {\rm (\ref{ppp1})}), and this
solution is classical, i.e. belongs to the Sobolev space
$H^{2}_{2}(M)$ if and only if the condition
$$
\frac{p(p-1)}{r_p^{2}}+\frac{q(q-1)}{r_q^{2}}>\frac{2(n-1)}{n-2}
$$
is fulfilled. 
\item[{\bf (2)}] If $\alpha\leq\frac{n+4}{4}<\alpha^*$, the equation
{\rm (\ref{cone1-39})} does not have other radial solutions in $L_2(K)$.
\item[{\bf (3)}] For $Q_{\alpha}>0$, $\alpha>\frac{n+4}{4}$, there
exists a family of solutions $C_F$ (Fowler solutions) which belong to
the space $L_2(M)$, but not to $H_2^1(M)$.
\end{enumerate}
\end{Theorem}
We denote, as above,
$\displaystyle
Q_{\alpha}=\inf_{\phi\in H^1_2(M),\ \phi\neq 0} I_{\alpha}(\phi).
$
\begin{Theorem}\label{non2}
Let $M$ be a manifold with tame conical singularity as above.
\begin{enumerate}
\item[{\bf (1)}] If $1<\alpha< \alpha^*$, then a minimizing function
$u_{\alpha}(\ell)$ exists, belongs to the space $H_2^1(M)$ and to the
space $H^{2}_{2}(M)$ if the plus-case condition {\rm (\ref{plus})} is
met. Asymptotically, $u_{\alpha}(\ell)\sim \ell^{-\sigma}$ near $x_*$,
where $\sigma =\frac{n-2}{2}(\mu -1)$.
\item[{\bf (2)}] If $\alpha=\alpha^*$, then a minimizing function
$u_{\alpha}(\ell)$ (existing in $H_2^1(M)$) belongs to $H_2^1(M)$,
and, if the plus-case condition {\rm (\ref{plus})} is met, to the
space $H^{2}_{2}(M)$ and $u_{\alpha}(\ell)\sim \ell^{-\sigma}$ near
$x_*$.
\end{enumerate}
\end{Theorem}
{\bf Remark.}  Notice that the asymptotic behavior of the solution
$u_{s}(\ell )\in H^{1}_{2}(M)$ {\bf is the same for all} $\alpha$,
$1\leqq \alpha \leqq \alpha ^{*},$ including the linear case $\alpha
=1$ , compare Section \ref{s6}.
\vspace{2mm}

\noindent
{\bf \ref{s7}.7. Nonradial solutions.} Theorem \ref{non2}, describes
the asymptotic behavior of those radial solutions of (\ref{cone1-30})
which belong to the Sobolev space $H^{2}_{2}(M)$.  It is important to
determine whether general (non-radial) solutions have better or worse
asymptotics then that of $u_{s}(\ell )$.
\vspace{2mm}

Similar questions have been studied extensively for the solutions of
the Yamabe equation on $S^n$ with a finite number of singularities
(in our terms, when $p=n$, $\alpha =\alpha ^{*}$, see \cite{KMPS} for
the most recent results and references to earlier works). It is shown
in \cite{KMPS} that any nonradial solution asymptotically behaves as
a shift (by $t$) of uniquely defined radial solutions of the same
equation. There is clear evidence that the same holds in our case.
\vspace{2mm}

However, for our purposes we do not need such a result in full strength.
We restrict our attention to the case when a solution tends to zero as
$t\rightarrow +\infty$. Then a modification of the proof of
\cite[Ch. 4,Theorem 5]{Coppel} may be done, so that it will work in
our case.  That gives the following result:
\begin{Theorem} Let $u(\ell ,\theta ,\psi )$ be a solution of
$(51)_{\alpha }$ such that
$$\Vert
 u\Vert _{L_{2}(S^{p}\times S^{q}) }(\ell ) \rightarrow 0$$ as $\ell
\rightarrow 0$.  Then there exists a radial solution $u_{0}(\ell )$ of
(\ref{cone1-30}) with the
same property such that
$$
\Vert u(\ell ,\theta ,\psi )-u_{0}(\ell )\Vert =o(\Vert u_{0}(\ell
)\Vert )\ \ \ \ \ \mbox{as $\ell \rightarrow 0$}.
$$
\end{Theorem}
As a result, the nonradial solutions of (\ref{cone1-30}) that tend to
zero, approaching the singular point $x_*$ have, for all $\alpha >1$,
the same asymptotic behavior as corresponding radial solutions (which
also approach zero).
\vspace{2mm}

Recall that the original metric $g$ on the cone is given by $g_K =
d\ell^2 + \frac{r_p^2\ell^2}{2}d^2 \theta + \frac{r_q^2\ell^2}{2}d^2
\psi$. According to \cite{AB}, there exists a Yamabe metric $\check{g}$
on $M$. 
\begin{Corollary}\label{yamabe1}
Let $(M,g)$ be a manifold with tame conical singularity as above, and
let $\check{g}\in [g]$ be a Yamabe metric on $M$. Then
$$
\check{g}\sim \ell^{-\sigma} \(d\ell^2 + \frac{r_p^2\ell^2}{2}d^2
\theta + \frac{r_q^2\ell^2}{2}d^2\)
$$
near the singular point $x_*$.
\end{Corollary}
\section{Appendix}\label{ap}
Here we collect some necessary computations, most of which are very
simple.
\vspace{2mm}

\noindent
{\bf \ref{ap}.1. Function $\Lambda=\Lambda(p)$.} Recall that in the
appropriate coordinates the scalar curvature on the cone $C(S^p\times
S^q)$ is given as
$$
\begin{array}{c}
R_{g_K}(\theta,\psi,\ell)={\Lambda\over \ell^2}, \ \ \mbox{with} \ \
\Lambda=p(p-1)\frac{2-r_p^{2}}{r_p^{2}}+q(q-1)\frac{2-r_q^{2}}{r_q^{2}}-2pq.
\end{array}
$$
We let $n$ be fixed and study $\Lambda $ as function of $p$ only,
$\Lambda=\Lambda(p)$.
\vspace{2mm}

\noindent
Substituting $q=n-p-1$, we transform $\Lambda $ to the expression
$$
\begin{array}{rcl}
\Lambda (p)
&=& 
-(n-1)(n-2)+2\left[ \frac{p(p-1)}{r_p^{2}}+\frac{q(q-1)}{r_q^{2}}
\right] 
\\
\\
&=& 
 p(p-1)\frac{2-r_p^{2}}{r_p^{2}}+q(q-1)\frac{2-r_q^{2}}{r_q^{2}}-2pq
\\
\\
&=& 
2p^{2}\left( r_p^{-2}+r_q^{-2} \right) -2p\left[ \left( r_p^{-2}+r_q^{-2}
\right)+(n-2)\frac{2}{r_q^{2}}\right] 
\\
\\
&&\ \ \ \ \ \ \ \ \ \ \ \ \ \ \ \ \ \ \ \ \ \ \ \ \ \ \ \ \ \ \ \ \ 
+
(n-1)(n-2)\left(\frac{2}{r_q^{2}}-1\right).
\end{array}
$$
To get the first equality we split
$\frac{2-r_p^{2}}{r_p^{2}}=\frac{2}{r_p^{2}}-1$ and use the fact that
$p(p-1)+q(q-1)-2pq=(n-1)(n-2)$.
\vspace{2mm}

\noindent
{\bf \ref{ap}.2. Minimal value of $\Lambda (p)$.}
To find the minimum of function $\Lambda (p)$ for $1<p <n-1$ we calculate
$$
\frac{d \Lambda }{dp}=4p\left( r_p^{-2}+r_q^{-2} \right)-2\left[ \left(
r_p^{-2}+r_q^{-2} \right)+(n-2)\frac{2}{r_q^{2}}\right],
$$ 
and find 
$$
p_{\min}=\frac{\left( r_p^{-2}+r_q^{-2}
\right)+(n-2)\frac{2}{r+q^{2}}}{2\left( r_p^{-2}+r_q^{-2}
\right)}=\frac{1}{2}+\frac{n-2}{1+\frac{r_q^{2}}{r_p^{2}}}.
$$
Correspondingly, $q$ takes the minimal value
$$
q_{\min}=\frac{1}{2}+\frac{n-2}{1+\frac{r_p^{2}}{r_q^{2}}}.
$$
Notice that for all possible values of $r_p>0$, $r_q>0$, we have
$$
\frac{1}{2} <p_{\min} <n-\frac{3}{2}.
$$
For the value of $\Lambda (p_{\min})$ we have
$$
\begin{array}{rcl}
\Lambda (p_{\min}) &=&
-(n-1)(n-2)+\frac{2}{r_p^{2}}\left(
\frac{1}{2}+\frac{n-2}{1+\frac{r_q^{2}}{r_p^{2}}}\right)
\left( \frac{-1}{2}+\frac{n-2}{1+\frac{r_q^{2}}{r_p^{2}}}\right) 
\\
\\
& &
\ \ \ \ \ \ \ \ \ \ \ \ \ \ \ \ \ \ \ \ \ \ + \ 
\frac{2}{r_q^{2}}\left( \frac{1}{2}+
\frac{n-2}{1+\frac{r_p^{2}}{r_q^{2}}}\right)
\left( \frac{-1}{2}+\frac{n-2}{1+\frac{r_p^{2}}{r_q^{2}}}\right) 
\\
\\
&= &
-(n-1)(n-2)-\frac{1}{2}( r_p^{-2}+r_q^{-2})+
\frac{2(n-2)^{2}}{r_p^{2}+r_q^{2}}.
\end{array}
$$
In Section \ref{s7}.5 we use the parameter $\mu ^{2}$.  We have:
$$
\begin{array}{rcl}
\mu ^{2}&=&\displaystyle  1+\frac{\Lambda }{(n-1)(n-2)}=
\\
\\
&= &\displaystyle 
1+\frac{p(p-1)\frac{2-r_p^{2}}{r_p^{2}}+q(q-1)
\frac{2-r_q^{2}}{r_q^{2}}-2pq }{(n-1)(n-2)}=
\\
\\
&= &\displaystyle 
\frac{ p(p-1)\frac{2}{r_p^{2}}+q(q-1)\frac{2}{r_q^{2}} }{(n-1)(n-2)} > 0
\end{array}
$$
if $p, q>0$.
Here we have used $p(p-1)+q(q-1)+2pq=(n-1)(n-2).$
\vspace{2mm}

\noindent
The maximal value is achieved at the end point of the interval
$[3,n-3]$ of the admissible values of $p$. We have for $p=1$,
$$
\begin{array}{rcl}
\Lambda(1) &=&\displaystyle  -(n-1)(n-2)+(n-2)(n-3)\frac{2}{r_q^{2}}, 
\\
\\
\Lambda(n-2)  &=& \displaystyle 
-(n-1)(n-2)+(n-2)(n-3)\frac{2}{r_p^{2}}.
\end{array}
$$
Correspondingly, maximal value of $\Lambda (p)$ is achieved at one of
these ends.
\vspace{2mm}

\noindent
{\bf \ref{ap}.2. Parameters $\bar{a}$, $\bar{b}$ in (\ref{cone1-42}).}
In order to analyze the ``phase portraits'' of the system
(\ref{cone1-42}) for different values of $Q_{\alpha}$ and $\Lambda$,
we need particular information on the parameters $\bar{a}$, $\bar{b}$.
Introduce the parameter $s=\frac{(\alpha-1)(n-2)}{4}$. Then it is easy
to see that $0< s \leq 1$ since $1<\alpha
\leq\alpha^*=\frac{n+2}{n-2}$.
\begin{Claim}\label{app1} For $n\geq 3$, $\bar{b}< 0$.
\end{Claim}
Indeed, for $n\geq 3$, we have
$$
\bar{b}= (n-2)-\frac{4}{\alpha-1} = (n-2)\(1-\frac{4}{(\alpha-1)(n-2)}\)=
- \frac{(n-2)(s-1)}{s} <0 .
$$
Now we study dependence of $\bar{a}$ on $s$ and $\Lambda$. We have
$$
\begin{array}{rcl}
\bar{a} &=& \displaystyle
\frac{(n-2)}{4(n-1)}\Lambda + \frac{2(n-2)}{\alpha-1} - \frac{4}{(\alpha-1)^2}
\\
\\
&=& \displaystyle
\frac{(n-2)}{4}\cdot\frac{\Lambda}{n-1} + \frac{4(n-2)^2}{2(\alpha-1)(n-2)} -
\frac{4\cdot 4(n-2)^2}{4(\alpha-1)^2(n-2)^2}
\\
\\
&=& \displaystyle
\frac{(n-2)}{4}\cdot\frac{\Lambda}{n-1} + \frac{(n-2)^2}{2s} - \frac{(n-2)^2}{4s^2}
\\
\\
&=& \displaystyle
\frac{(n-2)^2}{4}\[\frac{\Lambda}{(n-1)(n-2)}+ \frac{2}{s}-\frac{1}{s^2}\].
\end{array}
$$
To analyze the equilibrium point $(0,0)$, we need to know the sign of
expression $\frac{\bar{b}^2}{4}+\bar{a}$. We have (using expressions
for $\bar b$ an $\bar a$ obtained above):
$$
\begin{array}{rcl}
\displaystyle
\frac{\bar{b}^2}{4}+\bar{a} &=& \displaystyle
\frac{(n-2)^{2}(1-s)^{2}}{4s^{2}} +
\frac{(n-2)^2}{4}\[\frac{\Lambda}{(n-1)(n-2)}+ \frac{2}{s}-\frac{1}{s^2}\]
\\
\\
 &=& \displaystyle
\frac{(n-2)^{2}}{4}\left[ 1+\frac{\Lambda }{(n-1)(n-2)}\right]
=
\frac{(n-2)^{2}}{4}\mu ^{2}
\geq   0.
\end{array}
$$
Thus, we have that $\bar{a}(s) \to -\infty$ as $s \to 0$, and
$$
\bar{a}(1) = \frac{(n-2)^2}{4}\[\frac{\Lambda}{(n-1)(n-2)}+ 1\] =
\frac{(n-2)^{2}}{4}\mu^{2}> 0
$$
for all $\Lambda$.  Also, we have
$$
\frac{d\bar{a}}{ds} = \frac{(n-2)^2}{4}\[ -\frac{2}{s^2}+\frac{2}{s^3}\] =
\frac{(n-2)^2}{2s^3}\[1-s\]
$$
which is positive for all $s<1$, and zero for $s=1$. 
\vspace{2mm}

\parbox{1.8in}{\PSbox{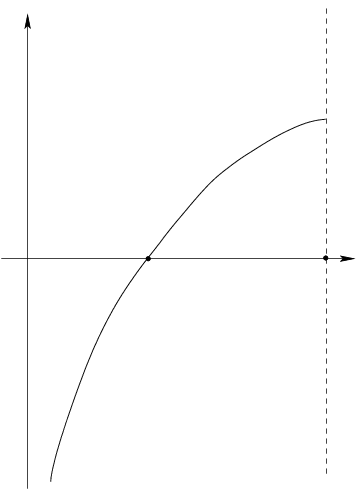}{8cm}{50mm}
\begin{picture}(0,0)
\put(10,150){{\small $\bar{a}(s)$}}
\put(40,73){{\small $s_0$}}
\put(87,73){{\small $1$}}
\put(105,83){{\small $s$}}
\end{picture}
\centerline{{\small {\bf Fig. \ref{ap}.1.}}}}
\parbox{2.75in}{Thus a graph
of the function $\bar{a}(s)$ has a form given at Fig. \ref{ap}.1. The
point $s_0$ (when $\bar{a}(s_0)=0$ is a root of the equation
$$
\frac{\Lambda}{(n-1)(n-2)}s^2 + 2s -1 = 0.
$$
One has the following positive root:
$$
s_0 = \frac{1}{1+\sqrt{1+\frac{\Lambda}{(n-1)(n-2)}}}=\frac{1}{1+\mu}
$$
for all $\Lambda$.}
\vspace{2mm}

\noindent
In terms of the parameter $\alpha$, $\bar{a}$ changes
sign from negative to positive if $\alpha=\alpha_0$, with
\begin{equation}\label{app2-q3}
\alpha_0 = 1 +
\frac{4}{(n-2)\(1+\sqrt{1+\frac{\Lambda}{(n-1)(n-2)}}\)}=1+\frac{4}{(n-2)(1+\mu)}.
\end{equation}
\begin{Claim}\label{app2}
The parameter $\bar{a}<0$ if $\alpha\in (1,\alpha_0)$, $\bar{a}>0$ if
$\alpha\in (\alpha_0,\alpha^*)$, and $\bar{a}=0$ if $\alpha=\alpha_0$, where
$\alpha_0$ is given by {\rm (\ref{app2-q3}).}
\end{Claim}
For the expression $\frac{{\bar b}^2}{4}-{\bar a}(\alpha -1)$ we get
$$
\begin{array}{l}
\frac{{\bar b}^2}{4}-{\bar a}(\alpha
-1) = 
\frac{(n-2)^{2}(s-1)^{2}}{4s^{2}}-\frac{4s}{(n-2)}\frac{(n-2)^{2}}{4}\left[
\frac{\Lambda }{(n-1)(n-2)} +\frac{2}{s}-\frac{1}{s^{2}}\right] 
\\
\\
= \! 
\left[\!  \frac{(n-2)^{2}}{4}-2(n-2) \right]\!  +\! 
\frac{(n-2)^{2}}{4}\frac{1}{s^{2}}\! +\! 
\left[ (n-2)\! -\! \frac{(n-2)^{2}}{2}\right] \frac{1}{s}\! 
-\! (n-2)(\mu ^{2}-1)s.
\end{array}
$$
As $s \rightarrow 0$, the leading term is the first one, and the value
of the function goes to $+\infty $. On the other hand, at $s=1$, the
value of this function is equal to $-(n-2)\mu ^{2} <0 $. Thus, this
expression takes both negative and positive values.
\vspace{2mm}

It was shown, in Lemma \ref{new-l1}, that a function with the $\alpha
$-basic asymptotic behavior at $x_*$ belongs to $L_{2}(M)$ iff
$\frac{n+4}{n} < \alpha $. It is instructive to compare this condition
with the condition on the $sign({\bar a})$. 
\vspace{2mm}

The condition
$\frac{n+4}{n}>\alpha _{0}$ is satisfied if and only if 
$\frac{n-2}{n} >
\frac{1}{1+\mu }$, i. e. iff
$$
\mu > \frac{2}{n-2}.
$$ 
Taking the square of this inequality and using the expression for $\mu $
we see that this condition is equivalent to
$$
\Lambda
>(n-1)(n-2)\left[\frac{4}{(n-2)^{2}-1}\right]=-\frac{n(n-1)(n-4)}{n-2}.
$$
Now we analyze the necessary and sufficient condition from Proposition
\ref{new-p1} when a solution $u$ belongs to the Sobolev space
$H^{k}_{2}$. The condition is
$$
k< 1\pm \sqrt{\frac{{\bar b}^{2}}{4}+{\bar a}}=1\pm \frac{n-2}{2}\mu .
$$
Here we reformulate this condition in terms of $p,P,q,Q$.
Consider first the minus-solutions, for which the condition takes the form
$$
k< 1-\frac{n-2}{2}\mu .
$$
Since $\mu >0$, such solution can not belong to $H^{1}_{2}$. On the
other hand, it belongs to $L_{2}$ if and only if $\frac{n-2}{2}\mu
<1$, or when $\mu <\frac{2}{n-2}$.  Substituting this expression for
$\mu $, we get the result
\begin{Claim}
The solution $u$ belongs to $L_{2}$ iff
$$
\frac{p(p-1)}{r_p^{2}}+\frac{q(q-1)}{r_q^{2}}<\frac{2(n-1)}{n-2}.
$$
\end{Claim}
{\bf Example:} For $r_p=r_q=1$ this condition takes the form
$(n-1-p)^{2}<2(n-1)$, or $p<n-1-\sqrt{2(n-1)}$.
\vspace{2mm}

\noindent
Consider now plus-solutions where the condition is $k<
1+\frac{n-2}{2}\mu.$ It is clear ($\mu >0$) that the solution $u$ always
belongs to $H^{1}_{2}$. This solution belongs to $H^{2}_{2}$ (i.e. is
the classical solution) iff $\mu <\frac{2}{n-2},$ or, substituting
the expression for $\mu $, when
$$
\frac{p(p-1)}{r_p^{2}}+\frac{q(q-1)}{r_p^{2}}>\frac{2(n-1)}{n-2}.
$$ 
\begin{Lemma}\label{app-new1}
If ${\bar a} <0$, a function with the $\alpha $-basic asymptotic
behavior at $x_*$ does not belong to $L_{2}(M).$
\end{Lemma}

\vspace{10mm}

\begin{small}
{\sf
\noindent
Boris Botvinnik, University of Oregon, Eugene
\\
e-mail: 
botvinn@math.uoregon.edu
\\
\\
Serge Preston, Portland State University,
\\
serge@mth.pdx.edu
}
\end{small}
\end{document}